%% file: lp-final-version.tex
\documentclass[11pt,reqno,a4]{amsart}

\usepackage[mathscr]{eucal}
\usepackage{enumerate}
\usepackage{epsfig,latexsym,amsfonts,amssymb,amsmath,amscd,
graphics, color, 
epic}
\newcommand{\bp}{\boldsymbol{p}}\newcommand{\bz}{\boldsymbol{z}}
\newcommand{\ba}{\boldsymbol{a}}\newcommand{\bm}{\boldsymbol{m}}\newcommand{\bn}{\boldsymbol{n}}
\newcommand{\by}{\boldsymbol{y}}\newcommand{\bu}{\boldsymbol{u}}\newcommand{\bv}{\boldsymbol{v}}
\newcommand{\bx}{\boldsymbol{x}}\newcommand{\bw}{\boldsymbol{w}}
\newcommand{\K}{\ensuremath{\mathbb{K}}}\newcommand{\bq}{\boldsymbol{q}}\newcommand{\ld}{\operatorname{ld}}
\newcommand{\Z}{\ensuremath{\mathbb{Z}}}
\newcommand{\Q}{\ensuremath{\mathbb{Q}}}\newcommand{\Tor}{\operatorname{Tor}}
\newcommand{\Aff}{\operatorname{Aff}}\newcommand{\Aut}{\operatorname{Aut}}
\newcommand{\Hessian}{\operatorname{Hessian}}
\newcommand{\R}{\ensuremath{\mathbb{R}}}
\newcommand{\C}{\ensuremath{\mathbb{C}}}
\newcommand{\tr}{\text{\rm tr}}
\newcommand{\pr}{\operatorname{pr}}
\newcommand{\PP}{\ensuremath{\mathbb{P}}}\newcommand{\Ini}{\operatorname{Ini}}
\newcommand{\Graph}{\operatorname{Graph}}\newcommand{\Trop}{\operatorname{Trop}}

\newcommand{\Vol}
{\ensuremath{\mathrm{Vol}}}

\newcommand{\mt}{\ensuremath{\mathrm{mt}}}

\newcommand{\conv}{\operatorname{Conv}}
\newcommand{\Sing}{\ensuremath{\mathrm{Sing}}}\newcommand{\eps}{\ensuremath{\varepsilon}}
\newcommand{\Val}{\ensuremath{\mathrm{Val}}}
\newcommand{\Span}{\ensuremath{\mathrm{Span}}}

\newcommand{\proofend}{\hfill$\Box$\bigskip}

\newtheorem{lemma}{Lemma}[section]
\newtheorem{theorem}[lemma]{Theorem}
\newtheorem{problem}[lemma]{Problem}
\newtheorem{corollary}[lemma]{Corollary}
\theoremstyle{remark}
\newtheorem{remark}[lemma]{Remark}

\newtheorem{example}[lemma]{Example}

\begin{document}
\title[Enumeration of Complex and Real Surfaces via Tropical Geometry]{Enumeration of Complex and Real Surfaces via Tropical Geometry}
\author {Hannah Markwig}
\address {Hannah Markwig, Eberhard Karls Universit\"at T\"ubingen, Fachbereich Mathematik, Institut f\"ur Geometrie, Auf der Morgenstelle 10, 72076 T\"ubingen, Germany }
\email {hannah@math.uni-tuebingen.de}

\author{Thomas Markwig}
\address{Thomas Markwig, Eberhard Karls Universit\"at T\"ubingen, Fachbereich Mathematik, Auf der Morgenstelle 10, 72076 T\"ubingen, Germany}
\email {keilen@math.uni-tuebingen.de}

\author{Eugenii Shustin}
\address{Eugenii Shustin, School of Mathematical Sciences, Tel Aviv University, Ramat Aviv, Tel Aviv 69978, Israel}
\email {shustin@post.tau.ac.il}

\begin{abstract} We prove a correspondence theorem for
singular tropical surfaces in $\R^3$, which recovers
singular algebraic surfaces in an appropriate toric three-fold that
tropicalize to a given singular tropical surface. Furthermore,
we develop a three-dimensional version of Mikhalkin's lattice path algorithm
that enumerates singular tropical surfaces passing through
an appropriate configuration of points in $\R^3$. As application
we show that there are pencils of real surfaces of degree $d$ in $\PP^3$
containing at least $(3/2)d^3+O(d^2)$ singular surfaces, which is
asymptotically comparable to the number $4(d-1)^3$ of all complex singular surfaces in the pencil.
Our
result relies on the classification of singular tropical surfaces
\cite{MMS}.
\end{abstract}
\keywords{Tropical geometry, tropical singular surface, discriminants, real singular surfaces}
\subjclass[2010]{14T05, 51M20, 14N10}
\maketitle

\tableofcontents

\section{Introduction}
\subsection{Main goal}
The tropical approach to enumerative geometry, initiated by Mikhalkin's correspondence theorem \cite{Mi2}, has led to remarkable success in the study of
Gromov-Witten and Welschinger (open Gromov-Witten) invariants of toric varieties (see, for example, \cite{IKS,Mi2}).
Mikhalkin originally used tropical methods to count curves in toric surfaces satisfying point conditions. Nowadays the tropical
techniques are developed for enumeration of
curves satisfying tangency conditions in addition \cite{BBM, GM}, for covers satisfying ramification conditions \cite{CJM} and for
curves in higher dimensional varieties satisfying point conditions \cite{NS}. Little is known about the enumerative geometry of surfaces in
toric three-folds and the tropical counterparts.
With this paper, we contribute a first step towards the establishment of tropical methods in such higher-dimensional enumerative problems.

The goal of this paper is to extend the tropical technique
to the case of surfaces in toric three-folds both in the complex and real setting, having in mind the test problem of
enumeration of complex and real surfaces belonging to a given
divisor class in
a given projective toric threefold, having a singularity in the big torus, and passing through an
appropriate number of generic points.
Even in this simply looking problem,
the tropical enumerative geometry appears to be non-trivial in each step.

Throughout the text, we say that a \emph{tropical surface is singular} if it is the tropicalization of a singular surface. A point on a singular tropical surface is called a \emph{tropical singular point}, or just singular point, if it is the tropicalization of a singular point. Notice that a tropical surface $S$ can have multiple singularities arising as tropicalizations of singular points of several surfaces in the fiber of $S$ of tropicalization.

\subsection{Main results}
The first result established in this paper is 
a three-dimensional version of the lattice path algorithm, which enumerates singular tropical surfaces

with a given Newton polytope that pass through a collection of points arranged on a generic line in a special way (see Lemmas
\ref{lBCE}, \ref{lD1}, \ref{lD2}, \ref{lD3}, \ref{lA1}, \ref{lA2},
\ref{lA3}, and \ref{lA4} in Section \ref{sec-circuit}).
Its idea is very similar to the original Mikhalkin's lattice path algorithm for plane tropical curves
(cf.
\cite{Mi2}). However, in the case of tropical surfaces,
the corresponding lattice paths can be disconnected, and the problem of inscribing
one of the five possible circuits dual to the face with a singular point of a singular
tropical surface\footnote{Recall that a \emph{circuit} is a set of
     lattice points that is affinely dependent but such that each
     proper subset is affinely independent.} (see \cite{MMS}) into the given lattice path turns to be a non-trivial combinatorial task,
which results both in local (i.e., related to the circuit) and global (i.e., related to the whole lattice path)
restrictions. Notice also that if the line of point constraints is sufficiently close to one of the
coordinate axes, in the planar curve case, the lattice path algorithm converges to
a Caporaso-Harris algorithm \cite{GM} which counts tropical
curves with relatively simple circuits represented by unit parallelograms and
multiple edges, and this leads
to a much simpler floor-diagram algorithm \cite{BM}. In its turn,
a similar enumeration of singular tropical surfaces necessarily involves surfaces
(see \cite[Appendix]{MMS2})
with circuits represented by unit parallelograms and double edges as well as by
pentatopes, which are much more involved (for example,
their classification up to $\Aff(\Z^3)$-equivalence is infinite, cf. \cite{MMS}).

The next result is a correspondence statement, which to a given configuration of points $\bp$ in the big torus of
the given toric three-fold and a singular tropical surface $S$ passing through the tropicalized point configuration
$\bx$
associates all possible singular algebraic surfaces which contain the configuration $\bp$ and
tropicalize to $S$, in particular, we compute the number of such algebraic surfaces, called
the \emph{multiplicity} $\mt(S,\bx)$ of the pair $(S,\bx)$. Thus, we obtain the formula
for the degree of the discriminant $\Sing(\Delta)$ in the tautological linear system on the toric variety associated with a convex lattice polytope $\Delta$:
\begin{equation}\deg\Sing(\Delta)=\sum_{S\in\Sing^\tr(\Delta,\overline\bx)}\mt(S,\overline\bx)\ .
\label{ecor03}\end{equation}

The combinatorial side of the correspondence
is addressed in Section \ref{sing-trop}, and it
consists
in the finding of all possible locations of tropical singular points on a given
singular tropical surface (see Lemmas \ref{lmt4} and \ref{lmt5} in Section
\ref{sing-trop}). Contrary to the planar case (see \cite{MMS,MMS1}),
these locations form a non-trivial finite set (cf. \cite{MMS}). As a byproduct we prove that
the twice lattice diameter of the Newton polytope provides
a lower bound to the maximal number of singular points on a singular tropical surface
with a given Newton polytope (Theorem \ref{t-lower} in Section \ref{sec-lower}).
The algebraic part of the correspondence is a version of Viro's patchworking construction
\cite{Vir}
adapted to the singular setting. That is we look for a singular surface over $\K$ as an
analytic family of complex surfaces inscribed in the family of three-folds ${\mathfrak X}\to(\C,0)$
with general fibre ${\mathfrak X}_t\simeq\Tor_\C(\Delta)$, $t\ne0$.
It is done in two steps: (1) in Section \ref{enhanced},
using the given configuration $\overline\bp\subset(\K^*)^3$, we restore the
the leading terms
of the coefficients of the polynomial describing the sought surface, i.e.,
the complex surface ${\mathcal S}_0\subset{\mathfrak X}_0$; (2)
in Section \ref{sec-pw}, we reconstruct the required family using the implicit function theorem.

In fact, there are even simpler formulas for the degree of the discriminant than (\ref{ecor03})
(see, for instance, \cite[Corollary
6.5]{DFS}). However, an advantage of the tropical approach is that it allows one to explicitly
describe all singular algebraic surfaces in count and, moreover, to recognize
the real ones among them. 
To demonstrate this advantage in our situation, we address the
following

\smallskip

{\bf Question:} {\it How many real singular surfaces can occur in a generic real pencil of surfaces of
degree $d$ in $\PP^3$?}

\smallskip

We show (Theorem \ref{t-real}
in Section \ref{sec-examples}) that there exist generic pencils of surfaces of degree $d$ containing at least
$(3/2)d^3+O(d^2)$ real singular surfaces, which is comparable with the total number
$4(d-1)^3$ of
(complex) singular surfaces in a pencil. We obtain also a general lower bound for the case of arbitrary Newton polytopes
(Theorem \ref{t-real1} in Section \ref{sec-examples}).


\medskip

\subsection{Acknowledgements}
The research was supported by the German-Israeli Foundation grant no. 1174-197.6/2011,
by the Minerva-Minkowski Center for Geometry at the Tel
Aviv University and by the DFG-grant MA 4797/5-1. A substantial part of this work has been
done during the authors' visit to the Centre Interfacultaire Bernoulli, Lausanne, and during the
third author's visit to Institut des Hautes \'Etudes Scientifiques. We are very grateful to
these institutions for
the hospitality and excellent working conditions. The authors would like to thank the unknown referee
for valuable remarks and suggestions, which helped us to improve the presentation.

\section{Preliminaries}\label{sec2}

\subsection{General setting} 
In the paper, we address complex and real enumerative problems, which can be expressed as counting intersections
of the discriminant with a complex or real pencil of hypersurfaces. Since the discriminant of a polynomial is defined
over $\Z$, by transfer principles (Lefschetz principle for the algebraically closed case
\cite[Theorem 1.13]{JL} and Tarski principle for the real closed case \cite[Theorem 1.16]{JL}), our problems
are respectively equivalent to those over
$\K=\bigcup_{m\ge1}\C\{t^{1/m}\}$, the field of locally convergent Puiseux series, and its real part
$\K_\R=\bigcup_{m\ge1}\R\{t^{1/m}\}$. Observe that $\K$ and $\K_\R$ possess the
non-Archimedean valuation
$\Val\left(\sum_ra_rt^r\right)=-\min\{r\in\Q\ :\ a_r\ne0\}\ .$

Let $\Delta\subset\R^3$ be a 
convex lattice polytope such that the set
$\{\bu-\bu'\ :\ \bu,\bu'\in\Delta\cap\Z^3\}$
generates the lattice $\Z^3$. Let $N=|\Delta\cap\Z^3|-2>0$. Denote by
$\Tor_\K(\Delta)$ the toric variety
over $\K$
associated to the polytope $\Delta$. Let ${\mathcal L}_\Delta$ be the
tautological line bundle on $\Tor_\K(\Delta)$.
Sections of ${\mathcal L}_\Delta$ are (Laurent) polynomials with support inside $\Delta$.
Denote by $|{\mathcal L}_\Delta|$ the linear system of divisors of non-zero elements
$\varphi\in H^0({\mathcal L}_\Delta)$. Clearly, $\dim|{\mathcal L}_\Delta|=|\Delta\cap\Z^3|-1=N+1$.
Define the discriminant $\Sing(\Delta)\subset|{\mathcal L}_\Delta|$ to be
the family parameterizing  divisors with a singularity in $(\K^*)^3$.
Assume that $\Delta$ is non-defective, i.e., the discriminant $\Sing(\Delta)$ is a hypersurface.
It is then natural to ask for its degree $\deg\Sing(\Delta)$.
The answer is known and can be expressed in terms of the combinatorics of $\Delta$ (see, for example,
\cite[Corollary 6.5]{DFS}). Sometimes it takes a form of a simple formula,
for instance, if
$\Delta$ is the simplex with vertices $(0,0,0)$, $(d,0,0)$, $(0,d,0)$ and$(0,0,d)$,
then
\begin{equation}\deg\Sing(\Delta)=4(d-1)^3\ .\label{ex1}\end{equation}

Geometrically, the degree can be seen as $\#(\Sing(\Delta)\cap{\mathcal P})$, where ${\mathcal P}\subset
|{\mathcal L}_\Delta|$ is a generic pencil. For example, we can take the pencil
${\mathcal P}_{\overline\bp}=\{{\mathcal S}\in|{\mathcal L}_\Delta|\ :\ {\mathcal S}
\supset\overline\bp\}\ ,$ where
$\overline\bp=(\bp_1,...,\bp_N)$ is a configuration of $N$ points in
$(\K^*)^3$ in general position. We intend to
explicitly describe the set $\Sing(\Delta)\cap{\mathcal P}_{\overline\bp}$ and, particularly, apply this
description to enumeration of real singular surfaces.


Denote by $\Sing^\tr(\Delta)$ the tropical discriminant parameterizing
singular tropical surfaces with
Newton polytope $\Delta$, i.e.\ tropicalizations of  algebraic surfaces ${\mathcal S}\in\Sing(\Delta)$
(background on singular tropical hypersurfaces can be found in \cite{DFS,DT,MMS}). Suppose that $\overline\bx=(\bx_1,...,
\bx_N) 
\subset\Q^3$ is a configuration of $N$ distinct points, which are in
general position, that is, the set $\Sing^\tr
(\Delta,\overline\bx):=\{S\in\Sing^\tr(\Delta)\ :\ S\supset\overline\bx\}$ is finite,
all tropical surfaces in this set are of maximal-dimensional geometric type (i.e.,
whose deformation space has dimension $\#(\Delta\cap\Z^3)-2$,
which is the maximal possible value for singular surfaces with Newton polytope $\Delta$), and the points $\bx_1,...,\bx_N$ are
interior points of $2$-faces in each of these surfaces
(see \cite[Theorem 1 and Section 2.3]{MMS}). 
Then we suppose that $\overline\bp\subset(\K^*)^3$ satisfies $\Val(\overline\bp)=\bx$ is generic among the configurations tropicalizing to $\overline\bx$.
In what follows we solve the following concrete problems:

\begin{problem}\label{main}
(1) Assuming that the configuration $\overline\bx\subset\Q^3$ is in Mikhalkin's position
(see Section \ref{sec-Mikhpos}), describe the combinatorics of tropical surfaces $S\in\Sing^\tr(\Delta)$ passing through
$\overline\bx$. We denote the set of these tropical surfaces by $\Sing^\tr(\Delta,\overline\bx)$.

(2) Given a tropical surface $S\in\Sing^\tr(\Delta,\overline\bx)$, calculate
$\mt(S,\overline\bx)$, the cardinality of the set
$\Sing(\Delta,\overline\bp,S)$ of surfaces ${\mathcal S}\in\Sing(\Delta)$ that
tropicalize to $S$ and pass through a fixed generic
configuration $\overline\bp\subset(\K^*)^3$ of $N$ points such that $\Val(\overline\bp)=\overline\bx$.

(3) Furthermore, assuming that the configuration $\overline\bp$ is real, calculate
$\mt^\R(S,\overline\bx)$, the number of real surfaces in $\Sing(\Delta,\overline\bp,S)$.
\end{problem}

\subsection{Classification of singular tropical surfaces}
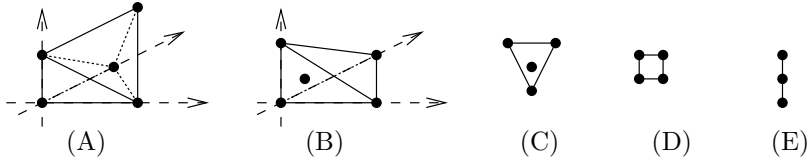
\begin{figure}[h]
    \centering
    \input{3circ.pstex_t}
    \caption{The possible circuits in the dual subdivision of a singular tropical surface.}
    \label{fig:circuits1}
  \end{figure}

Our results rely on the classification of singular tropical surfaces in $\R^3$
of maximal-dimensional geometric type \cite[Theorem 2]{MMS}. For the reader's convenience, we
present here this
classification:

\begin{lemma}\label{lrev1}
The dual subdivision of a singular tropical surface $S$ of a maximal-dimensional geometric type has a unique
circuit $C_S$, and the cell in the tropical surface dual to $C_S$ contains the
tropical singular point. Possible circuits are as follows (cf. Figure \ref{fig:circuits1}):
\begin{enumerate}\item[({\bf A})]
a pentatope whose vertices are its only integer points, equivalent up to $\Aut(\Z^3)$-action to one of the following
$$\{(0,0,0),\ (1,0,0),\ (0,1,0),\ (0,0,1),\ (1,p,q)\},\quad 1\le p\le q,\ (p,q)=1\ ,$$
\item[({\bf B})] a tetrahedron whose integer points are the vertices and
one interior point, equivalent up to $\Aut(\Z^3)$-action to one of the following (asterisk marks the interior point)
$$\{(0,0,0),(1,0,0),(0,1,0),(1,1,1)^*,(3,3,4)\},\quad \{(0,0,0),(1,0,0),(0,1,0),(1,1,2)^*,(2,2,5)\}\ ,$$
$$\{(0,0,0),(1,0,0),(0,1,0),(1,2,3)^*,(2,4,7)\},\quad \{(0,0,0),(1,0,0),(0,1,0),(1,3,5)^*,(2,6,11)\}\ ,$$
$$\{(0,0,0),(1,0,0),(0,1,0),(1,3,5)^*,(2,7,13)\},\quad \{(0,0,0),(1,0,0),(0,1,0),(1,4,7)^*,(2,9,17)\}\ ,$$
$$\{(0,0,0),(1,0,0),(0,1,0),(1,5,7)^*,(2,13,19)\},\quad \{(0,0,0),(1,0,0),(0,1,0),(1,2,5)^*,(3,7,20)\}\ ,$$
\item[({\bf C})] a triangle whose integer points are the vertices and one interior point, equivalent up to
$\Aut(\Z^3)$-action to $\{(1,0,0),(0,1,0),(1,1,0),(2,2,0)\}$,
\item[({\bf D})] a parallelogram whose vertices are the only its integral points, equivalent up to
$\Aut(\Z^3)$-action to $\{(0,0,0),(1,0,0),(0,1,0),(1,1,0)\}$,
\item[({\bf E})] a segment of lattice length $2$, equivalent up to $\Aut(\Z^3)$-action to
$\{(0,0,0),(1,0,0),(2,0,0)\}$.\end{enumerate}
\end{lemma}

\section{The lattice path algorithm in dimension $3$}\label{sec-lp}
In this section, we present a solution for Problem \ref{main} (1), which consists in the following algorithmic procedure:
\begin{itemize}\item First, we choose tropical point constraints in Mikhalkin's position
(see Section \ref{sec-Mikhpos}).
\item Next, we enumerate all possible lattice paths of length $N=|\Delta\cap\Z^3|-2$
inscribed into the polytope $\Delta$ and related to the chosen point constraints
(see Lemma \ref{llp1} in Section \ref{subsec:dualreformulation}).
\item Finally, for each of the above lattice paths and each of the five types of
circuits (Lemma \ref{lrev1}), we construct all singular tropical surfaces
that pass through the given point constraints, have a circuit of the chosen type, and
whose dual subdivision of $\Delta$ contains the given lattice path (see Lemmas
\ref{lBCE}, \ref{lD1}, \ref{lD2}, \ref{lD3}, \ref{lA1}, \ref{lA2},
\ref{lA3}, and \ref{lA4} in Section \ref{sec-circuit}).
\end{itemize}

In what follows we use the notation $\Vol_\Z(\delta)$ for the lattice volume of
a positive-dimensional lattice polytope $\delta$, i.e., the volume normalized by the condition that
the minimal lattice simplex of dimension $\dim\delta$ in the affine space
spanned by $\delta$ has volume $1$.

\subsection{Tropical point constraints in Mikhalkin's position}\label{sec-Mikhpos}
%

To apply a lattice path algorithm similar to the one for tropical curves \cite{Mi1}, \cite[Section 7.2]{Mi2}, we place the points in the
following special
position. Choose a line $L\subset\R^3$ passing
through the origin
and directed by a vector $\bv\in\Q^3$, which is not parallel or orthogonal to any proper affine
subspace of $\R^3$ spanned by a non-empty subset $A\subset\Delta\cap\Z^3$;
then pick the following (ordered) configuration $\overline\bx=(\bx_1,...,\bx_N)$ of marked points
\begin{equation}\begin{cases}&\bx_i=
M_i\bv\in L,\quad i=1,...,N,\quad\text{where}\\
&0\ll M_1\ll...\ll M_N\
\text{are positive rationals},
\end{cases}\label{eCH}\end{equation}
$N=|\Delta\cap\Z^3|-2$.

\begin{remark}\label{rCH} (1) In what follows in Section
\ref{sec-lp}, we consider $M_1,...,M_N$ in (\ref{eCH}) as
parameters to be chosen follows. There will appear finitely many linear combinations of
$M_1,...,M_N$ with coefficients depending only on the polytope $\Delta$, and we always assume that
if $\lambda(M_1,...,M_k)$ is such a combination for some $1\le k\le N$, where the coefficient of $M_k$
is positive, then $\lambda(M_1,...,M_k)\gg0$. We also use the notation
$\lambda=\Theta(M_k)$ if $\lambda=\lambda(M_1,...,M_k)$ is a linear combination as above
with a non-zero coefficient of $M_k$, and we write $\lambda=
o(M_k)$ if $\lambda$ is a linear combination of $M_i$'s with $i<k$. Thus, our assumption yields that
always $|o(M_k)|\ll|\Theta(M_k)|$.

(2) Observe that the configurations (\ref{eCH}) are generic. The set
$\Sing^\tr(\Delta,\overline\bx)$
is finite, and all its elements are singular tropical surfaces of maximal-dimensional geometric type
as described in \cite[Theorem 2]{MMS}.
Moreover, for any $S\in\Sing^\tr(\Delta,\overline\bx)$,
each marked point $\bx_i$, $1\le i\le N$ is in the interior of a $2$-face $F_i$ of $S$,
and $F_i\ne F_j$ as $i\ne j$.
\end{remark}


We will solve Problem \ref{main}(1) for point configurations satisfying (\ref{eCH}).

\subsection{The dual reformulation}\label{subsec:dualreformulation} Introduce the partial order in $\R^3$:
$\bu\succ\bu'\: \Longleftrightarrow\ \langle\bu-\bu',\bv\rangle>0$ to obtain a
linear order on $\Delta\cap\Z^3$:
$$\Delta\cap\Z^3=\{\bw_0,...,\bw_{N+1}\},\quad \bw_i\prec \bw_{i+1}\ \text{for all}\ i=0,...,N\ .$$
Given a subset $A\subset\Delta\cap\Z^3$, consisting of $m\ge2$ points
$\ba_1\prec\ba_2\prec...\prec\ba_m$, we call an ordered subset of the set of segments
$P(A):=\{[\ba_i,\ba_{i+1}]\ :\ i=1,...,m-1\}$ a \emph{lattice path} supported on $A$ if it covers the
whole set $A$. The set $P(A)$ is called the \emph{complete lattice path} supported on $A$.
We call a lattice path connected (disconnected) if the union of its segments is connected
(disconnected). 

Let $F_S:\R^3\to\R$ be a
tropical polynomial defining a singular tropical surface $S \in \Sing^\tr(\Delta,\overline\bx)$, $\nu_S:\Delta\to\R$ the Legendre dual
piecewise linear function, whose linearity domains determine
the subdivision $\Sigma_S$ of $\Delta$ dual to $S$.
Denote by $e_i$, $i=1,\ldots,N$ the edge of $\Sigma_S$ dual to the $2$-face $F_i$ of $S$ containing the point $\bx_i$ in its interior.
We denote by \mbox{$P(S,\overline\bx)=\{e_i\ :\ i=1,...,N\}$} the lattice path corresponding to the pair
$(S,\overline\bx)$.

\begin{lemma}\label{llp1}
 For a singular tropical surface $S$ passing through $\overline \bx$, the lattice path $P(S,\overline\bx)$ defined above satisfies:
 \begin{enumerate}
  \item[(i)] Either $P(S,\overline\bx)=P(A')\cup P(A'')$, where $A'=\{\bw_0,...,\bw_k\}$,
  $A''=\{\bw_{k+1},...,\bw_{N+1}\}$ for some $1\le k\le N$; we call this path $\Gamma_{k,k+1}$;
  \item[(ii)] or $P(S,\overline\bx)=P(A)$, where $A=\Delta\cap\Z^3\setminus\{\bw_k\}$ for some
  $0\le k\le N+1$; we call this path $\Gamma_k$.
 \end{enumerate}
\end{lemma}

We call the lattice paths $\Gamma_k$, $k=0,\ldots,N+1$ and $\Gamma_{k,k+1}$, $k=1,\ldots,N$ the \emph{marked lattice paths} for $\Delta$.

{\bf Proof.}
By the duality of $S$ and the subdivision $\Sigma_S$ (see \cite[Section 2.1]{Mi3}), the components of
$\R^3\setminus S$ are in one-to-one correspondence with a subset
of $\Delta\cap\Z^3$ (including
all the vertices of $\Delta$). Due to the convexity of these components,
different connected components of $L\setminus\overline\bx$ cannot intersect the same component of $\R^3\setminus S$. Since $L\setminus\overline\bx$ has
$|\overline\bx|+1=N+1=|\Delta\cap\Z^3|-1$ components, we encounter
the following situations:
\begin{enumerate}
 \item[(a)] both $L\setminus\overline\bx$ and $\R^3\setminus S$ consist of
 $N+1$ components;
 \item[(b)] $L\setminus\overline\bx$ consists of $N+1$ components, and
 $\R^3\setminus S$ consists of $N+2$ components.
\end{enumerate} Now note that if
$\bw_i$ and $\bw_j$ are dual to the components $\bw_i^*$, $\bw_j^*$
intersecting $L$ along neighboring intervals, and the vector $\bv$ points from
$\bw_i^*$ to $\bw_j^*$, then $\bw_j\succ\bw_i$.

In case (a), there exists a unique point $\bw_k$, $0\le k\le N+1$,
that is not a vertex of the subdivision $\Sigma_S$. Then $P(S,\overline\bx)=\Gamma_k$.

In case (b), if there is a component $\bw_k^*$ of $\R^3\setminus S$ disjoint from
$L$, then $P(S,\overline\bx)$ again is $\Gamma_k$ for some $0\le k\le N+1$. Otherwise,
we have an extra intersection point $\by\in L\setminus\overline\bx$ of $L\cap S$, and then:
if $\by\prec\bx_1$ we get the path $\Gamma_0$,
if $\bx_k\prec\by\prec\bx_{k+1}$ for some $k=1,...,N$, we get the path $\Gamma_{k,k+1}$, and
at last, if $\bx_N\prec\by$ we get the path
$\Gamma_{N+1}$.
\proofend

We can now refine problem \ref{main}(1) as follows:

\begin{problem}\label{pr1}
Given a marked lattice path $P$, find all subdivisions $\Sigma$ of $\Delta$ that contain the path $P$
(i.e., each edge of $P$ is an edge of the subdivision $\Sigma$) and are dual to
singular tropical surfaces $S$ passing through $\overline \bx$ (such that the edge dual to the $2$-face $F_i$ containing $\bx_i$ is in $P$).
\end{problem}

%
%

We suggest a solution to Problem \ref{pr1}, which can be regarded as a three-dimensional version of
Mikhalkin's lattice path algorithm \cite{Mi1,Mi2}. By \cite[Theorem 2]{MMS}, the desired subdivision $\Sigma$
has one circuit of type {\bf A}, {\bf B}, {\bf C}, {\bf D}, or {\bf E} (see
Lmma \ref{lrev1} and Figure \ref{fig:circuits1}) and all its three-dimensional cells that do not contain the circuit are
simplices, i.e.\ tetrahedra whose only integral points are their vertices. In the next Section \ref{sec-smooth},
we present an auxiliary construction that completes the subdivision outside the circuit.
In Section \ref{sec-circuit}, we explain how to fit a circuit in a subdivision
for a given lattice path.

\subsection{The smooth extension algorithm}\label{sec-smooth}

We will first show in general terms, how to extend a given subdivision
when the underlying polytope is enlarged.

\begin{lemma}\label{smooth-dop}
Let us be given the following data:
\begin{itemize}
\item a convex lattice polytope $\delta'\subset\R^n$
and a convex piecewise linear function $\nu':\delta'\to\R$, whose linearity domains
define a subdivision $\sigma'$ of $\delta'$
into convex lattice subpolytopes;
\item a convex lattice polytope $\delta''\subset\R^n$ such that $\delta_0=\delta'\cap\delta''$ is
a cell of the subdivision $\sigma'$ and a face of $\delta''$ of codimension $1$.
\end{itemize}
Pick a point $\bw\in\delta''\cap\Z^n\setminus\delta'$. Then there exits a unique extension
of $\sigma'$ to a convex subdivision $\sigma$ of $\delta=\conv(\delta'\cup\delta'')$ such that
\begin{itemize}\item the
vertices of $\sigma$ are the vertices of $\sigma'$ and of $\delta''$,
\item $\delta''$ is a cell of $\sigma$, \item
the cells of $\sigma$ are linearity domains of a convex piecewise linear function $\nu:\delta\to\R$ such that
$\nu\big|_{\delta'}=\nu'$ and $\nu(\bw)\gg\max\nu'$.\end{itemize}
\end{lemma}

{\bf Proof.}
Clearly, $\bw$ does not lie in the affine subspace of $\R^n$ spanned by $\delta_0$. Hence the (linear)
function $\nu'\big|_{\delta_0}$ and the value $\nu(\bw)$ induce a unique linear function $\nu''$ on
$\delta''$. Furthermore, the condition $\nu(\bw)\gg\max\nu'$ ensures that any segment in $\R^{n+1}$ joining
an interior point of the graph of $\nu'$ and an interior point of the graph of $\nu''$ lies above these graphs.
Hence the lower facets of $\conv(\Graph(\nu')\cup\Graph(\nu''))$ (i.e.\ the facets whose outer normal vector has a negative last coordinate) defines a graph of a convex piecewise linear function
$\nu:\delta\to\R$ as required. Finally, we note that there is a $\mu\gg\max\nu'$ such that the subdivision of
$\delta$ defined by the linearity domains of $\nu$ does not depend on the choice of the value $\nu(\bw)>\mu$.
\proofend

\begin{example}\label{ex-smooth}
Let $\delta'\subset\R^n$ and $\nu':\delta'\to\R$ be as in Lemma
\ref{smooth-dop}, $\bw\in\Z^n\setminus\delta'$, $\delta=\conv(\delta'\cup\{\bw\})$.
Let $\bv\in\Q^n$ be a vector which is not parallel or
orthogonal to any segment joining any two distinct points of $\delta$. Suppose that $\bw\succ\bw'$ for any
$\bw'\in\delta'$. Then the construction of Lemma \ref{smooth-dop} works as follows.
Note that there exists a point $\tilde{\bw}\in\delta'$ which satisfies $\tilde{\bw}\succ\bw'$ for all $\bw'\in\delta'
\setminus\{\tilde{\bw}\}$ and that the segment $[\tilde{\bw},\bw]$ intersects with $\delta'$ only at $\tilde{\bw}$. Then
we can put $\delta''=[\tilde{\bw},\bw]$ and extend the subdivision $\sigma'$ of $\delta'$ to a convex subdivision
of $\delta$. We call the subdivision $\sigma$ of $\delta$ the \emph{smooth extension of} $\sigma'$
\footnote{``Smooth" means here that the dual tropical surface cannot have singular points on the
faces dual to the cells of the subdivision lying outside $\delta'$.}.
\end{example}

An important particular case is the following construction.

\begin{lemma}\label{smooth} Let $\Delta=\conv(A)$, where $A\subset\Delta\cap\Z^3$, $|A|=N+1$, and
$A=\{\ba_0,...,\ba_N\}$,
$\ba_0\prec\ba_1\prec...\prec \ba_N$ (order defined by $\bv$).
Let $\overline\bx$ be a sequence of $N$ points of $\R^3$ given by (\ref{eCH}). Then:
\begin{enumerate}
\item[(i)] In the space of tropical surfaces defined by tropical polynomials of the
form
$$F:\R^3\to\R,\quad F(X)=\max_{\omega\in A}(c_i+\langle\ba_i,X\rangle),\quad c_i\in\R,\ i=0,...,N\ ,$$
there exists a unique surface $S=S(A,\overline\bx)$, that passes through
$\overline\bx$.
\item[(ii)] Each point of $\overline\bx$ belongs to the interior of some $2$-face of $S$, and distinct
points belong to distinct faces.
\item[(iii)] The dual subdivision $\Sigma_S$ consists of only tetrahedra, and it is constructed
by a sequence of smooth extensions, when starting with the point $\ba_0$ and subsequently adding
the points $\ba_1,...,\ba_N$. The edges dual to
the faces of $S$, that intersect $\overline\bx$, form the lattice path $P(A)$ subsequently going through
the points $\ba_0,...,\ba_N$.
\end{enumerate}
\end{lemma}

Notice that we view the space of tropical surfaces defined by tropical polynomials as above as $\R^{|\mathcal{A}|}/(1,\ldots,1)$. In particular, we can always assume that the first coefficient of the tropical polynomial satisfies $c_0=0$.

{\bf Proof.}
Statements (ii) immediately follows from the general position of $\overline\bx$.
Thus, we explain only parts (i) and
(iii).
The polynomial $F_S(X)$ defining $S$ can be computed from the formulas:
$$c_0=0,\quad c_{i-1}+\langle\ba_{i-1},\bx_i\rangle=c_i+\langle\ba_i,\bx_i\rangle,\ i=1,...,N\ ,$$
or, equivalently,
\begin{equation}c_0=0,\quad c_i-c_{i-1}=-M_i\langle\ba_i-\ba_{i-1},\bv\rangle,\ i=1,...,N\ .
\label{eA2}\end{equation} The function
$\nu_S:\Delta\to\R$ takes value $-c_i$ at the point $\ba_i$, $i=0,...,N$.
Since $0\ll M_1\ll...\ll M_N$, we have $\nu_S(\ba_i)\gg \nu_S(\ba_{i-1})$ for all $i=1,...,N$,
which is required in Lemma \ref{smooth-dop} and Example \ref{ex-smooth}.
\proofend

\subsection{Subdivisions with prescribed type of circuit}\label{sec-circuit}

In this section, we study how the types of circuits listed in
Lemma \ref{lrev1} fit into subdivisions for a given lattice path. In particular, we show that the circuits of type
{\bf B}, {\bf C}, {\bf D}, or {\bf E} can appear only for connected lattice paths $\Gamma_k$ omitting one
point $\bw_k\in\Delta\cap\Z^3$, which is either an interior points of the convex hull of the circuit of
type {\bf B}, {\bf C}, or {\bf E}, or is one of the points of the circuit of type {\bf D}. In turn,
circuits of type {\bf A} may appear in subdivisions based on connected and disconnected lattice paths.

\subsubsection{Subdivisions with circuit of type {\bf B}, {\bf C}, or {\bf E}}\label{secBCE}

\begin{lemma}\label{lBCE}
(1) A marked lattice path $P$ admits an extension to a subdivision $\Sigma$ of $\Delta$, dual to a surface
$S\in \Sing^\tr(\Delta,\overline\bx)$ and having a circuit of type {\bf B},
{\bf C}, or {\bf E}, only if $P=\Gamma_k$ (see Lemma \ref{llp1}), where
$1\le k\le N$, and $\bw_k$ is not a vertex of $\Delta$. Moreover, this subdivision
is unique and it can be constructed by the
smooth extension algorithm of Lemma \ref{smooth}(iii)
supported on the set $A=\Delta\cap\Z^3\setminus\{\bw_k\}$.

(2) Let $P=\Gamma_k$, where $1\le k\le N$ and $\bw_k$ is not a vertex of $\Delta$. Then the subdivision
$\Sigma$ of $\Delta$, constructed as in item (1), is dual to a surface $S\in\Sing^\tr(\Delta,\overline\bx)$ if and
only if one of the following conditions holds true:
\begin{itemize}
\item the point $\bw_k$ belongs to the interior of a three-dimensional cell of $\Sigma$
(i.e.\ $\bw_k$ is the interior point of a circuit of type {\bf B});
\item the point $\bw_k$ belongs to the interior of a two-dimensional cell of $\Sigma$, and,
if $\bw_k\in\partial\Delta$, the subdivision $\Sigma$ additionally satisfies the third condition
in \cite[Theorem 4]{MMS} (i.e.\ $\bw_k$ is the interior point of a circuit of type {\bf C});
\item the point $\bw_k$ is the midpoint of an edge of $\Sigma$, and, if $\bw_k\in\partial\Delta$, the
subdivision $\Sigma$ additionally satisfies the fourth condition in \cite[Theorem 4]{MMS}
(i.e.\ $\bw_k$ is the interior point of a circuit of type {\bf E}).
\end{itemize}
\end{lemma}

{\bf Proof.}
Statement (1) is straightforward.
Statement (2) follows from \cite[Theorem 4]{MMS}.
\proofend

\begin{remark}\label{rem-circuitrelation}
It follows from the smooth extension algorithm of Lemma \ref{smooth}(iii) that the coefficients $c_i$, $i\neq k$, of the tropical polynomial defining the unique surface $S\in \Sing^\tr(\Delta,\overline\bx)$ dual to a subdivision extending $\Gamma_k$ and containing a circuit of type {\bf B}, {\bf C}, or
{\bf E} according to Lemma \ref{lBCE}(2) are determined by the point conditions $\overline \bx$. Furthermore, the lattice points $\bw_l$ forming the circuit satisfy a unique up to nonzero multiple relation $\sum_l \lambda_l \bw_l=0$ with $\sum \lambda_l=0$
(for example, for circuit of type {\bf E} it is
$\bw_i-2\bw_j+\bw_k=0$ with $\bw_j$ the midpoint of the segment $[\bw_i,\bw_k]$). Since the circuit is part of the subdivision, it follows that $\sum_l \lambda_l c_l=0$, which allows us to deduce the value of $c_k$ from the others. We call the equation $\sum_l \lambda_l c_l=0$ defining $c_k$ the \emph{circuit relation} for the coefficients of the tropical polynomial.
\end{remark}

\subsubsection{Subdivisions with circuit of type {\bf D}}\label{secD}

For circuits of type {\bf D}, we have to treat the case of a connected path $\Gamma_k$ or a disconnected path $\Gamma_{k,k+1}$ (see Lemma \ref{llp1}) separately.

{\bf (1) The case of a connected path $P$.}

\begin{lemma}\label{lD1}
Let $P=\Gamma_k$ for some $k=0,...,N+1$, and let $P$ extend to a subdivision
$\Sigma$ of $\Delta$ with a circuit $C$ of type {\bf D}, that is dual to a surface
$S\in\Sing^\tr(\Delta,\overline\bx)$. Then
\begin{enumerate}
\item[(i)] the circuit $C$ contains $\bw_k$ and three more vertices
$\bw_i,\bw_j,\bw_l$, $i<j<l$;
\item[(ii)] the subdivision $\Sigma$ is uniquely determined by the pair $(k,C)$, in particular,
\begin{itemize}
 \item it contains a smooth triangulation
of $\conv(P(l^*))$ as in Lemma \ref{smooth}, where $P(l^*)$ is the part of $P$ bounded from above by the
vertex $\bw_{l^*}$ preceding $\bw_l$ in $P$,
\item the parallelogram $\conv(C)$ intersects $\conv(P(l^*))$ along the edge $[\bw_i,\bw_j]$,
\item $\Sigma$ is obtained from the triangulation of $\conv(P(l^*))$ by the
extension to $\conv(P(l^*)\cup C)$ as in Lemma \ref{smooth-dop} and by a sequence of smooth
extensions as in Example \ref{ex-smooth} when subsequently adding the points of $P$ following $\bw_l$.
\end{itemize}
\end{enumerate}
\end{lemma}

{\bf Proof.} We explain only the first claim in statement (ii), since the rest is straightforward.
As in Lemma \ref{smooth} and Remark \ref{rem-circuitrelation}, we obtain the coefficients of a tropical polynomial defining the surface $S$ from the point conditions and the circuit relation. Let $\nu_S$ be the piece-wise linear function defined by this polynomial.

Suppose that $\bw_s\in P$, $\Delta_s=\conv(P(s))$ (where $P(s)$ is the part of $P$ bounded from above by the
vertex $\bw_s$) is smoothly triangulated, and $\bw_s$ is the maximal such lattice point. Assume that $s<l^*$. The fact that the triangulation of $\Delta_s$ does not extend to a smooth triangulation
of $\conv(\Delta_s\cup\{\bw_{s+1}\})$ means that in the graph of
$\nu$, there exists a line segment $\sigma_1$ joining $(\bw_{s+1}, \nu_S(\bw_{s+1}))$
with a point $(\bz_1,\nu_S(\bz_1))\in\Delta_s\times \R$ and a line segment $\sigma_2$ joining
a point $(\bw_m,\nu_S(\bw_m))$, $m>s+1$, or the point $(\bw_k,\nu_S(\bw_k))$ with a point $(\bz_2,\nu_S(\bz_2))\in\Delta_s\times \R$,
such that $\sigma_1\cap(\Delta_s\times \R)=(\bz_1,\nu_S(\bz_1))$, $\sigma_2\cap(\Delta_s\times \R)=(\bz_2,\nu_S(\bz_2))$,
$\sigma_2$ lies in a lower face of the graph of $\nu_S$, and the
projections of $\sigma_1,\sigma_2$ onto $\R^3$ intersect
in the interior of the projection of this face. This, however, contradicts the convexity of the
function $\nu_S:\Delta\to\R$, since the values $\nu_S(\bw_m)$, where $m>s+1$ or $m=k$,
are much larger than $\nu_S(\bw_{s+1})$ (for $m\neq k$ this follows from the smooth extension algorithm Lemma \ref{smooth}, for $m=k$ from the circuit relation as in Remark \ref{rem-circuitrelation}).
\proofend

Lemma \ref{lD1} provides only necessary conditions for a connected lattice path with circuit of type
{\bf D}
to be extendable to a subdivision, dual to a singular
tropical surface passing through the given point configuration. To formulate sufficient conditions,
consider the univariate tropical polynomial
\begin{equation}F_S\big|_L(\tau)=\max_{0\le s\le N+1}(c_s+\tau\langle\bw_s,\bv\rangle)\ .
\label{eA3}\end{equation}
Its coefficients $c_0,...,c_{N+1}$ are
determined by the following relations (point conditions (\ref{eCH}) and circuit relation, see Remark \ref{rem-circuitrelation}):
\begin{itemize}
\item for $k=0$
$$c_1=0,\quad c_{s+1}+M_s\langle\bw_{s+1},\bv\rangle=c_s+M_s\langle\bw_s,\bv\rangle,
\ 1\le s\le N,\quad c_0+c_l=c_i+c_j\ ,$$
\item for $k=N+1$
$$c_0=0,\quad c_s+M_s\langle\bw_s,\bv\rangle=c_{s-1}+M_s\langle\bw_{s-1},\bv\rangle,
\ 1\le s\le N,\quad c_i+c_{N+1}=c_j+c_l\ ,$$
\item for $1\le k\le N$
$$c_0=0,\quad\begin{cases}c_s+M_s\langle\bw_s,\bv\rangle=c_{s-1}+M_s\langle\bw_{s-1},\bv\rangle,\ &\text{as}\ 1\le s<k,\\
c_{k+1}+M_k\langle\bw_{k+1},\bv\rangle=c_{k-1}+M_k\langle\bw_{k-1},\bv\rangle,\ &\\
c_{s+1}+M_s\langle\bw_{s+1},\bv\rangle=c_s+M_s\langle\bw_s,\bv\rangle,\ &\text{as}\ k<s\le N,\end{cases}$$
$$\begin{cases}c_k+c_l=c_i+c_j,\ &\text{if}\ k<i,\\
c_i+c_l=c_k+c_j,\ &\text{if}\ i<k<l,\\
c_i+c_k=c_j+c_l,\ &\text{if}\ k>l.\end{cases}$$
\end{itemize}

\begin{lemma}\label{lD2}
The subdivision $\Sigma$ constructed in Lemma \ref{lD1} is dual to a tropical surface $S\in\Sing^\tr(\Delta,\overline\bx)$
if and only if the following conditions hold:
\begin{enumerate}
\item[(i)] the face of $\Sigma$ given by the circuit does not lie on $\partial\Delta$;
\item[(ii)] $i<k<l$.
\end{enumerate}
\end{lemma}

{\bf Proof.}
The first condition is necessary by \cite[Theorem 4]{MMS}. Having it fulfilled, we have to ensure that
the roots of the tropical polynomial $F_S\big|_L(\tau)$ (that is the restriction of the
tropical polynomial defining $S$ to the line $L$) are $M_s$, $1\le s\le N$, and
maybe one more root outside the range $[M_1,M_N]$. Since the tropical polynomial
$$\widetilde F(\tau)=\max_{0\le s\le N+1,\ s\ne k}(c_s+\tau\langle
\bw_s,\bv\rangle)$$ has precisely the roots $M_1,...,M_N$, we
end up with inequalities \begin{equation}\begin{cases}c_0+M_1\langle\bw_0,\bv\rangle\le M_1\langle\bw_1,\bv\rangle
=c_2+M_1\langle\bw_2,\bv\rangle,
\quad &\text{if}\ k=0,\\ c_k+M_k\langle\bw_k,\bv\rangle\le c_{k-1}+M_k\langle\bw_{k-1},\bv\rangle
\\
\qquad\qquad\qquad\qquad=c_{k+1}+M_k\langle\bw_{k+1},\bv\rangle,
\quad &\text{if}\ 1\le k\le N,\\
c_{N+1}+M_N\langle\bw_{N+1},\bv\rangle\le c_N+M_N\langle\bw_N,\bv\rangle
\\ \qquad\qquad\qquad\qquad\qquad=c_{N-1}+M_N\langle\bw_{N-1},\bv\rangle,
\quad &\text{if}\ k=N+1.\end{cases}\label{eD2}\end{equation}

Condition (\ref{eD2}) is necessary as well,
cf. the proof of Lemma \ref{llp1}. The sufficiency of conditions (i) and (\ref{eD2}) comes again from
\cite[Theorem 4]{MMS} (i.e., $S$ is singular), and from the fact that
$F_S\big|_L(\tau)$ defines the intersection points of $L$ and $S$:

We show the equivalence of (ii) and (\ref{eD2}).
Suppose that $i<j<l<k\le N$. Then from (\ref{eD2}) and the circuit relation we derive
\begin{equation}c_i=o(M_l),\ c_j=o(M_l),\ c_l=\Theta(M_l),
\ c_{k-1}=o(M_k),\ c_k=c_l+c_j-c_i=\Theta(M_l),\label{eD3}\end{equation} and hence
$$c_k+M_k\langle\bw_k,\bv\rangle\le c_{k-1}+M_k\langle\bw_{k-1},\bv\rangle\quad
\Leftrightarrow\quad c_k\le c_{k-1}-M_k\langle\bw_k-\bw_{k-1},\bv\rangle$$
$$\Longrightarrow\quad \Theta(M_l)\le-M_k\langle\bw_k-\bw_{k-1},\bv\rangle+o(M_k)\ ,$$ a
contradiction.

Suppose that $i<j<l<N<k=N+1$. Again (\ref{eD2}) and the circuit relations yield
$$c_{N+1}+M_N\langle\bw_{N+1},\bv\rangle\le c_N+M_N\langle\bw_N,\bv\rangle\quad\Longleftrightarrow
\quad c_{N+1}\le c_N-M_N\langle\bw_{N+1}-\bw_N,\bv\rangle$$
$$\Longrightarrow\quad \Theta(M_l)\le-M_N\langle\bw_{N+1}-\bw_{N-1},\bv\rangle+o(M_N)\ ,$$
a contradiction, since $l<N$, and hence
$c_l=\Theta(M_l)=o(M_N)$.

Suppose that $i<j<l=N<k=N+1$. Then similarly we get
$$c_{N+1}+M_N\langle\bw_{N+1},\bv\rangle\le c_N+M_N\langle\bw_N,\bv\rangle\quad\Longleftrightarrow
\quad c_{N+1}\le c_N-M_N\langle\bw_{N+1}-\bw_N,\bv\rangle$$
$$\Longleftrightarrow\quad
c_N+c_j-c_i\le c_N-M_N\langle\bw_{N+1}-\bw_N,\bv\rangle
\quad\Longleftrightarrow\quad c_j-c_i\le-M_N\langle\bw_{N+1}-\bw_N,\bv\rangle$$
which is a contradiction.

Suppose that $1\le k<i<j<l$. Then we have
\begin{equation}\begin{cases} &c_i=o(M_{l-1}),\ c_j=o(M_{l-1}),\ c_l=
-M_{l-1}\langle\bw_l-\bw_{l-1},\bv\rangle+o(M_{l-1}),\\
& c_{k+1}=o(M_{l-1}),\ c_k=c_i+c_j-c_l=
M_{l-1}\langle\bw_l-\bw_{l-1},\bv\rangle+o(M_{l-1}),\end{cases}
\label{eD4}\end{equation}
and hence
$$c_k+M_k\langle \bw_k,\bv\rangle\le c_{k+1}+M_k\langle\bw_{k+1},\bv\rangle\quad
\Longleftrightarrow\quad c_k\le c_{k+1}+M_k\langle\bw_{k+1}-\bw_k,\bv\rangle$$
$$\Longrightarrow\quad M_{l-1}\langle\bw_l-\bw_{l-1},\bv\rangle+o(M_{l-1})\le o(M_{l-1}),$$
a contradiction.

In the case $k=0<1<i<j<l$, we again have relations
(\ref{eD4}), and hence
$$c_0+M_1\langle\bw_0,\bv\rangle\le M_1\langle\bw_1,\bv\rangle\quad
\Longleftrightarrow\quad c_0\le M_1\langle\bw_1-\bw_0,\bv\rangle$$
$$\Longrightarrow\quad M_{l-1}\langle\bw_l-\bw_{l-1},\bv\rangle+o(M_{l-1})\le o(M_{l-1}),$$
a contradiction.

In the case $k=0$, $i=1<j<l$, we similarly obtain
$$c_0\le M_1\langle\bw_1-\bw_0,\bv\rangle=\Theta(M_1)=o(M_{l-1})\ ,$$ since $l-1\ge2$. However,
from the circuit relation, we get
$$c_0=c_j-c_l=M_{l-1}\langle\bw_l-\bw_{l-1},\bv\rangle+o(M_{l-1})\ ,$$ which contradicts the former conclusion.

Suppose that $i<k<l$. Then the equations for $c_s$, $0\le s\le N+1$, yield
$$c_l=\begin{cases}-M_{l-1}\langle\bw_l-\bw_{l-1},\bv\rangle+o(M_{l-1}),\quad &\text{if}\ k<l-1,\\
-M_{l-1}\langle\bw_l-\bw_{l-2},\bv\rangle+o(M_{l-1}),\quad &\text{if}\ k=l-1.\end{cases}$$
If $k<l-1$, the required relation reads
$$c_k+M_k\langle\bw_k,\bv\rangle\le c_{k+1}+M_k\langle\bw_{k+1},\bv\rangle\quad
\Longleftrightarrow\quad c_k\le c_{k+1}+M_k\langle\bw_{k+1}-\bw_k,\bv\rangle=\Theta(M_k)=o(M_{l-1})$$
$$\Longleftrightarrow\quad c_l+c_i-c_j\le o(M_{l-1})\quad\Longleftrightarrow\quad
-M_{l-1}\langle\bw_l-\bw_{l-1},\bv\rangle+o(M_{l-1})\le o(M_{l-1})\ ,$$ which holds true.

If $k=l-1$, the required relation reads
$$c_l+c_i-c_j+M_{l-1}\langle\bw_{l-1},\bv\rangle\le c_l+M_{l-1}\langle\bw_l,\bv\rangle\quad
\Longleftrightarrow\quad c_i-c_j=o(M_{l-1})\le M_{l-1}\langle\bw_l-\bw_{l-1},
\bv\rangle\ ,$$ which again holds true.
\proofend

\smallskip

{\bf (2) The case of a disconnected path $P$.}

\begin{lemma}\label{lD3}
 Let $P=\Gamma_{k,k+1}$ for some
$1\le k<N$. Then it cannot be extended to a subdivision of
$\Delta$ with a circuit of type {\bf D} that is dual to a surface $S\in\Sing^\tr(\Delta,\overline\bx)$.
\end{lemma}

{\bf Proof.}
Observe that $L\cap S$ contains the marked points $\bx_s=M_s\bv$,
$1\le s\le N$, and one more point $\bx_0=M_0\bv$
such that $M_k<M_0<M_{k+1}$, which separates the intervals $\bw^*_k\cap L$ and $\bw^*_{k+1}\cap L$
(here, $\bw^*_k$ denotes the connected component of $\R^3\setminus S$ dual to $\bw_k$, cf. the proof of Lemma \ref{llp1}). If $P$ extends to a subdivision of $\Delta$ dual to a surface
$S\in\Sing^\tr(\Delta,\overline\bx)$, then the coefficients of
the tropical polynomial $F_S|_L(\tau)$
(see (\ref{eA3}))
can be computed from
\begin{equation}c_0=0,\quad c_{s+1}=\begin{cases} c_s-M_{s+1}\langle\bw_{s+1}-\bw_s,\bv\rangle,
                       \quad & \text{if}\ 0\le s<k,\\
                       c_k-M_0\langle\bw_{k+1}-\bw_k,\bv\rangle,\quad &
                       \text{if}\ s=k,\\
                       c_s-M_s\langle\bw_{s+1}-\bw_s,\bv\rangle,\quad &\text{if}
                       \ k<s\le N.
                      \end{cases}\label{eD5}\end{equation}
Assume the circuit of type {\bf D} consists of the points $\bw_i$, $\bw_j$, $\bw_l$ and $\bw_m$, with $i<j<l<m$. Joining relations (\ref{eD5}) and $0\ll M_1\ll...\ll M_N$, we see that 
the
circuit relation $c_i+c_m=c_j+c_l$ 
can hold only if $l=k$, $m=k+1$, and $M_0=\Theta(M_k)$. However, under these conditions, the
actual circuit relation $c_i+c_{k+1}=c_j+c_k$ converts to
$$c_{k+1}-c_k=c_j-c_i\quad\Longrightarrow\quad-M_0\langle\bw_{k+1}-\bw_k,\bv\rangle= \Theta(M_j)\ ,$$
which is a contradiction, since $M_0>M_k\gg M_j$ as $j<k$.
\proofend

\subsubsection{Subdivisions with circuit of type {\bf A}}\label{secA}
Recall (cf. \cite[Theorem 2]{MMS})
that a circuit of type {\bf A} is formed by the vertices of a pentatope, which up to
$\Z$-affine transformation can be identified with
\begin{equation}
\Pi_{p,q}=\conv\{(0,0,0),(1,0,0),(0,1,0),(0,0,1),(1,p,q)\},\quad p,q>0,\ \gcd(p,q)=1\ .
\label{epen}\end{equation}
The circuit relation (see Remark \ref{rem-circuitrelation}) means that the points $(\omega,-c_\omega)\in\R^4$, $\omega\in\Pi_{p,q}$,
$c_\omega\in\R$, lie in one $3$-plane, and it can be written as
\begin{equation}c_{100}+pc_{010}+qc_{001}=(p+q)c_{000}+c_{1pq}\ .\label{eA1}\end{equation}

\smallskip

{\bf (1) The case of a connected path $P$.}

\begin{lemma}\label{lA1}
 Let the lattice path $P=\Gamma_k$ (see Lemma \ref{llp1}), $0\le k\le N+1$, admit an extension to
 a subdivision $\Sigma$ of
 $\Delta$ with a circuit $C=\{\bw_i,\bw_j,\bw_l,\bw_m,\bw_n\}$, $i<j<l<m<n$,
 of type {\bf A} dual to a surface $S\in\Sing^\tr(\Delta,\overline\bx)$. Then
 \begin{enumerate}
  \item[(i)] $k\in\{i,j,l,m,n\}$;
  \item[(ii)] the cases $k=n\le N$ and $k=n=N+1>m+1$ are not possible;
  \item[(iii)] the subdivision $\Sigma$ is uniquely determined by the pair
  $(k,C)$ and satisfies the following:
  \begin{itemize}
   \item it contains a smooth
  triangulation of $\Delta_{m-1}=\conv\{\bw_s\ :\ 0\le s<m,\ s\ne k\}$;
  \item the pentatope $\conv(C)$ intersects $\Delta_{m-1}$ along their common
 $2$-face spanned by the first three points of $C\setminus\{\bw_k\}$;
  \item $\Sigma$ is obtained from the triangulation of $\Delta_{m-1}$ by the
extension to \mbox{$\conv(\Delta_{m-1}\cup C)$} as in Lemma \ref{smooth-dop} and by a sequence of smooth
extensions as in Example \ref{ex-smooth} when subsequently adding the points of $P$ following $\bw_n$.
  \end{itemize}
 \end{enumerate}
\end{lemma}

{\bf Proof.}
Claim (i) immediately follows from formulas (\ref{eA2}), since in
case $k\not\in\{i,j,l,m,n\}$, we would have
$|c_n|\gg\max\{|c_i|,|c_j|,|c_l|,|c_m|\}$ contrary to the circuit relation
(\ref{eA1}) (combined with a proper $\Z$-affine transformation).

Suppose now that $k=n\le N$. The necessary condition in this case is (see (\ref{eD2}))
$$c_n+M_n\langle\bw_n,\bv\rangle\le c_{n-1}+M_n\langle\bw_{n-1},\bv\rangle
$$ $$\Longrightarrow\quad
c_n\le c_{n-1}-M_n\langle\bw_n-\bw_{n-1},\bv\rangle=-M_n
\langle\bw_n-\bw_{n-1},\bv\rangle+o(M_n)\ ,$$
whereas from the circuit relation (\ref{eA1}) we get
$$c_n 
=\Theta(M_m)=o(M_n)\ ,$$ a contradiction.

Suppose that $k=n=N+1>m+1$. Then the necessary condition (\ref{eD2}) yields
$$c_{N+1}+M_N\langle\bw_{N+1},\bv\rangle\le c_{N-1}+M_N\langle\bw_{N-1},\bv\rangle$$
$$\Longrightarrow\quad c_{N+1}=c_{N-1}-M_N\langle\bw_{N+1}-\bw_{N-1},\bv\rangle=
-M_N\langle\bw_{N+1}-\bw_{N-1},\bv\rangle+o(M_N)\ ,$$
which again contradicts the circuit relation
$$c_{N+1} 
=\Theta(M_m)=o(M_N)\ .$$

Claim (iii) is proved analogously to Lemma \ref{lD1}(ii).
\proofend

\begin{lemma}\label{lA2}
In the notation of Lemma \ref{lA1}, let the data $k$ and $C$ satisfy conditions
(i) and (ii), and let a subdivision $\Sigma$ of $\Delta$ be constructed as in item (iii).
Write the circuit relation (\ref{eA1}) in the form
\begin{equation}c_k=\sum_{s\in\{i,j,l,m,n\}\setminus\{k\}}\lambda_sc_s\ .\label{eA4}\end{equation}
Then $\Sigma$ is dual to a tropical surface $S\in\Sing^\tr(\Delta,\overline\bx)$ if and only if
the following holds:
\begin{itemize}
 \item for $k=n=N+1$, $m=N$, either
 \begin{equation}\langle(\lambda_N-1)(\bw_N-\bw_{N-1})-(\bw_{N+1}-\bw_N),\bv\rangle>0\ ,\label{eA5}\end{equation}
 or
 \begin{equation}\begin{array}{l}\langle(\lambda_N-1)(\bw_N-\bw_{N-1})-(\bw_{N+1}-\bw_N),\bv\rangle=0
 \quad\text{and}\\ \qquad \begin{cases}
                      & \text{either}\  l<N-1,\\
                       & \text{or}\ l=N-1,\ \lambda_N-1+\lambda_{N-1}>0, \\
                       & \text{or}\ l=N-1,\ \lambda_N-1+\lambda_{N-1}=0,\ \lambda_j>0,
                      \end{cases}\end{array}
\label{eA7}\end{equation}
 \item for $0\le k<n$, we have $\lambda_n>0$.
\end{itemize}
\end{lemma}

{\bf Proof.}
Similarly to the proof of Lemma \ref{lD2}, we have to check conditions (\ref{eD2}), which read
\begin{equation}\begin{cases}c_{N+1}+M_N\langle\bw_{N+1},\bv\rangle\le c_N+M_N\langle
   \bw_N,\bv\rangle,\quad & \text{if}\ k=n=N+1,\ m=N,\\
   c_k+M_k\langle\bw_k,\bv\rangle\le c_{k-1}+M_k\langle\bw_{k-1},\bv\rangle,
   \quad &\text{if}\ 0<k<n,\\
   c_0+M_1\langle\bw_0,\bv\rangle\le M_1\langle\bw_1,\bv\rangle,\quad&\text{if}\
   k=0.
  \end{cases}\label{eA6}\end{equation}

  In the first case, we plug the circuit relation (\ref{eA4}) and the relation
  $c_N=c_{N-1}-M_N\langle\bw_N-\bw_{N-1},\bv\rangle$ into (\ref{eA6}) and obtain
  \begin{equation}(\lambda_N-1)c_{N-1}+\lambda_lc_l+\lambda_jc_j+\lambda_ic_i\le
  M_N\langle(\lambda_N-1)(\bw_N-\bw_{N-1})-(\bw_{N+1}-\bw_N),\bv\rangle\ . \label{eA8}\end{equation}
  Since the left-hand side is of order $o(M_N)$, we immediately see that
  (\ref{eA5}) is sufficient for (\ref{eA6}), and that the opposite strict inequality
  contradicts (\ref{eA6}). If the right-hand side of (\ref{eA8}) vanishes, we get $\lambda_N-1>0$,
  and hence conditions (\ref{eA7}) in view of
  $$c_{N-1}=-M_{N-1}\langle\bw_{N-1}-\bw_{N-2},\bv\rangle+o(M_{N-1}),\quad c_l=-M_l\langle
  \bw_l-\bw_{l-1},\bv\rangle+o(M_l),\quad c_j,c_i=o(M_l)\ .$$

In the second case, we again plug the circuit relation into (\ref{eA6}) and obtain
$$\begin{cases}\lambda_nc_n+\sum_{s\in\{i,j,l,m\}\setminus\{k\}}\lambda_sc_s-c_{k-1}\le
   -M_k\langle\bw_k-\bw_{k-1},\bv\rangle,\quad & \text{if}\ k\ne0,\\
   \lambda_nc_n+\lambda_mc_m+\lambda_lc_l+\lambda_jc_j\le
   M_1\langle\bw_1-\bw_0,\bv\rangle,\quad & \text{if}\ k=i=0,
  \end{cases}$$ which holds if and only if $\lambda_n>0$ in view of
  $$c_n=-M_n\langle\bw_n-\bw_{n-1},\bv\rangle+o(M_n),\quad c_i,c_j,c_l,c_m=o(M_n)$$
  (recall that $\lambda_n\ne0$).
\proofend

\smallskip

{\bf (2) The case of a disconnected path $P$.}

\begin{lemma}\label{lA3}
 Let the lattice path $P=\Gamma_{k,k+1}$, $1\le k\le N$, (see Lemma \ref{llp1}) admit an extension to a subdivision $\Sigma$ of
 $\Delta$ with a circuit $C=\{\bw_i,\bw_j,\bw_l,\bw_m,\bw_n\}$, $i<j<l<m<n$,
 of type {\bf A} dual to a surface $S\in\Sing^\tr(\Delta,\overline\bx)$. Then
 \begin{enumerate}
  \item[(i)] either $m=k$, $n=k+1$, or $m=k+1$, $n=k+2$;
  \item[(ii)] the subdivision $\Sigma$ is uniquely determined by the pair
  $(k,C)$ and satisfies the following:
  \begin{itemize}
   \item it contains a smooth
  triangulation of $\Delta_{m-1}=\conv\{\bw_s\ :\ 0\le s<m,\ s\ne k\}$;
  \item the pentatope $\conv(C)$ intersects $\Delta_{m-1}$ along their common
 $2$-face $\conv\{\bw_i,\bw_j,\bw_l\}$;
  \item $\Sigma$ is obtained from the triangulation of $\Delta_{m-1}$ by the
extension to \mbox{$\conv(\Delta_{m-1}\cup C)$} as in Lemma \ref{smooth-dop} and by a sequence of smooth
extensions as in Example \ref{ex-smooth} when subsequently adding the points of $P$ following $\bw_n$.
  \end{itemize}
 \end{enumerate}
\end{lemma}

{\bf Proof.}
From equations (\ref{eD5}), we get that
$c_i,c_j,c_l=o(|c_n|)$, and hence the circuit relation (\ref{eA1}) yields that $c_m$ and $c_n$ must be of the same
order. This is only possible if either $m=k$, $n=k+1$, and $M_0$ is comparable with $M_k$, or
$m=k+1$, $n=k+2$, and $M_0$ is comparable with $M_{k+1}$. Claim (ii) can be proved as Lemma \ref{lD1}(ii).
\proofend

\begin{lemma}\label{lA4}
In the notation of Lemma \ref{lA3}, let the data $k$ and $C$ satisfy condition
(i), and let a subdivision $\Sigma$ of $\Delta$ be constructed as in item (ii).
Write the circuit relation (\ref{eA1}) in the form
\begin{equation}c_n=\lambda_mc_m+\lambda_lc_l+\lambda_jc_j+\lambda_ic_i\ .\label{eA9}\end{equation}
Then $\Sigma$ is dual to a tropical surface $S\in\Sing^\tr(\Delta,\overline\bx)$ if and only if
the following holds:
\begin{itemize}
 \item for $m=k$, $n=k+1$, we have $\lambda_k-1>0$ and either
 \begin{equation}
 \langle\bw_{k+1}-\bw_k-(\lambda_k-1)(\bw_k-\bw_{k-1}),\bv\rangle<0\ ,
 \label{eA13}\end{equation}
 or \begin{equation}
 \begin{array}{l}
 \langle\bw_{k+1}-\bw_k-(\lambda_k-1)(\bw_k-\bw_{k-1}),\bv\rangle=0,\quad\text{and}\\
 \qquad\qquad\qquad\begin{cases}\text{either}\quad & l<k-1,\\
 \text{or}\quad & l=k-1,\ \lambda_k+\lambda_{k-1}-1>0,\\
 \text{or}\quad & l=k-1,\ \lambda_k+\lambda_{k-1}-1=0,\ \lambda_j>0,\end{cases}
 \end{array}\label{eA14}\end{equation}
 \item for $m=k+1$, $n=k+2$, we have $\lambda_{k+1}-1>0$ and either
 \begin{equation}\langle\bw_{k+2}-\bw_{k+1}-(\lambda_{k+1}-1)(\bw_{k+1}-\bw_k),\bv\rangle<0\ ,
 \label{eA17}\end{equation}
 or \begin{equation}\begin{array}{l}
 \langle\bw_{k+2}-\bw_{k+1}-(\lambda_{k+1}-1)(\bw_{k+1}-\bw_k),\bv\rangle=0,\quad\text{and}\\
 \qquad\qquad\begin{cases}\text{either}\quad & l<k,\\
 \text{or}\quad &l=k,\ \lambda_{k+1}+\lambda_k-1>0,\\
 \text{or}\quad &l=k,\ \lambda_{k+1}+\lambda_k-1=0,\ \lambda_j>0.\end{cases}
 \end{array}\label{eA18}\end{equation}
\end{itemize}
\end{lemma}

{\bf Proof.}
Suppose that $m=k$ and $n=k+1$. Plugging $c_{k+1}=c_k-M_0\langle\bw_{k+1}-\bw_k,\bv\rangle$ into
(\ref{eA9}), we obtain
$$-M_0\langle\bw_{k+1}-\bw_k,\bv\rangle=(\lambda_k-1)c_k+\lambda_lc_l+\lambda_jc_j+\lambda_ic_i\ .$$
The required condition $M_{k+1}>M_0>M_k$ is equivalent to the two inequalities:
\begin{itemize}\item $-M_{k+1}\langle\bw_{k+1}-\bw_k,\bv\rangle<(\lambda_k-1)c_k+\lambda_lc_l+\lambda_jc_j+\lambda_ic_i$,
which holds true, since by (\ref{eD5})
\begin{equation}c_k,c_l,c_j,c_i=o(M_{k+1})\ ;\label{eA19}\end{equation}
\item $-M_k\langle\bw_{k+1}-\bw_k,\bv\rangle>(\lambda_k-1)c_k+\lambda_lc_l+
\lambda_jc_j+\lambda_ic_i$, which via the substitution of
$c_k=c_{k-1}-M_k\langle\bw_k-\bw_{k-1},\bv\rangle$ transfers into
\begin{equation}
-M_k\langle\bw_{k+1}-\bw_k-(\lambda_k-1)(\bw_k-\bw_{k-1}),\bv\rangle>
(\lambda_k-1)c_{k-1}+\lambda_lc_l+\lambda_jc_j+\lambda_ic_i\ .
\label{eA12}\end{equation}
Since $c_{k-1},c_l,c_j,c_i=o(M_k)$ by (\ref{eD5}), we immediately get that $\lambda_k-1>0$, that (\ref{eA13})
is sufficient for (\ref{eA12}), and that
the opposite strict inequality in (\ref{eA13}) contradicts
(\ref{eA12}). At last, if the left-hand side of (\ref{eA12}) vanishes, due to
$$\begin{cases}&c_{k-1}=-M_{k-1}\langle\bw_{k-1}-\bw_{k-2},\bv\rangle+o(M_{k-1})<0,\\
&c_l=-M_l\langle\bw_l-\bw_{l-1},\bv\rangle+o(M_l)<0,\\
&c_j=-M_j\langle\bw_j-\bw_{j-1},\bv\rangle+o(M_j)<0,\end{cases}$$ we end up with condition (\ref{eA14}).
\end{itemize}

Suppose that $m=k+1$ and $n=k+2$. Plugging $c_{k+2}=c_{k+1}-M_{k+1}\langle\bw_{k+2}-\bw_{k+1},\bv\rangle$
and $c_{k+1}=c_k-M_0\langle\bw_{k+1}-\bw_k,\bv\rangle$ into
(\ref{eA9}), we obtain
$$(\lambda_{k+1}-1)M_0\langle\bw_{k+1}
-\bw_k,\bv\rangle-M_{k+1}\langle\bw_{k+2}-\bw_{k+1},
\bv\rangle=(\lambda_{k+1}-1)c_k+\lambda_lc_l+\lambda_jc_j+\lambda_ic_i\ .$$
Observe that this yields $\lambda_{k+1}-1>0$ in view of (\ref{eA19}).
Furthermore, we again have to satisfy the inequalities $M_{k+1}>M_0>M_k$, which are equivalent to:
\begin{itemize}\item
$(\lambda_{k+1}-1)M_k\langle\bw_{k+1}
-\bw_k,\bv\rangle-M_{k+1}\langle\bw_{k+2}-\bw_{k+1},
\bv\rangle<(\lambda_{k+1}-1)c_k+\lambda_lc_l+\lambda_jc_j+\lambda_ic_i$,
which always holds due to (\ref{eA19}); and
\item $M_{k+1}\langle(\lambda_{k+1}-1)(\bw_{k+1}-\bw_k)-
(\bw_{k+2}-\bw_{k+1}),\bv\rangle>(\lambda_{k+1}-1)c_k+\lambda_lc_l+\lambda_jc_j+\lambda_ic_i$, which
holds under condition (\ref{eA17}) and fails under the opposite strict inequality in
(\ref{eA17}) in view of (\ref{eA19}). Finally, if the left-hand side of (\ref{eA17}) vanishes, due to
$$\begin{cases}&c_k=-M_k\langle\bw_k-\bw_{k-1},\bv\rangle+o(M_k)<0,\\
&c_l=-M_l\langle\bw_l-\bw_{l-1},\bv\rangle+o(M_l)<0,\\
&c_j=-M_j\langle\bw_j-\bw_{j-1},\bv\rangle+o(M_j)<0,\end{cases}$$ we end up with condition (\ref{eA18}).\proofend
\end{itemize}


\section{Multiplicities of singular tropical surfaces}\label{sec-mt}

\subsection{General setting}\label{sec-setting}
Now, given a singular tropical surface $S\in\Sing^\tr(\Delta,\overline\bx)$,
we restore all singular algebraic surfaces over $\K$ with Newton polytope $\Delta$,
passing through the a generic configuration $\overline\bp\subset((\K^*)^3)^N$, $\Val(\overline\bp)=\overline\bx$, and
tropicalizing to $S$. In particular,
we compute their number $\mt(S,\overline\bx)$.
This number is finite due to the general position of the configuration $\overline\bp$,
but it may vanish as we see below, since the singular lifts
of a tropical surface $S\in\Sing^\tr(\Delta,\overline\bx)$ may avoid the configuration $\overline\bp$.

We follow the general patchworking procedure in the style of \cite[Chapter 2]{IMS} or
\cite{Sh}. It amounts to the following: (i) the tropical surface $S$ defines a toric degeneration of
the toric three-fold $\Tor_\C(\Delta)$, namely, a family ${\mathfrak X}\to(\C,0)$ with a general fiber
$\Tor_\C(\Delta)$ and the central fiber ${\mathfrak X}_0$ splitting into the union of toric three-folds determined by the subdivison of $\Delta$ dual to $S$; the point configuration
$\overline\bp$ (defined over $\K$) turns into the set of sections of the above family;
(ii) using the configuration $\overline\bp_0\subset{\mathfrak X}_0$ we find suitable (reducible)
algebraic surfaces ${\mathcal S}_0\subset{\mathfrak X}_0$ passing through $\overline\bp_0$;
(iii) finally, we extend each ${\mathcal S}_0$ to a family
${\mathcal S}\to(\C,0)$ of singular algebraic surfaces inscribed into the family ${\mathfrak X}\to
(\C,0)$ and containing the sections $\overline\bp$, i.e., we obtain
singular algebraic surfaces over the field $\K$ tropicalizing to $S$ and passing through $\overline\bp$.
Accordingly, we proceed in three steps:
\begin{enumerate}
\item[(1)] In Section \ref{sing-trop}, we find possible locations of singular points
in $S$; this is relevant for the case of
circuits of type {\bf C} and {\bf E}, for which the position of the tropical
singular points is not determined uniquely (Lemmas
\ref{lmt4} and \ref{lmt5}).
\item[(2)] In Section \ref{enhanced}, we describe the family ${\mathfrak X}\to(\C,0)$ and
find suitable surfaces ${\mathcal S}_0\subset{\mathfrak X}_0$ (Lemmas \ref{lmt2} and \ref{lmt3}). \item[(3)] In Section \ref{sec-pw}, we find the desired
singular algebraic surfaces in the form of families ${\mathcal S}\to(\C,0)$ (Lemmas
\ref{lmt11}, \ref{lmt6}, \ref{lmt8}, and \ref{lmt7}); notice that the data collected in steps (1) and (2)
do not determine the family ${\mathcal S}\to(\C,0)$ uniquely, so, we attach additional information (like
the position of the singular point in ${\mathcal S}_0$) which can be interpreted as an extra blowing up
of ${\mathfrak X}$ in order to obtain transversal conditions and finally apply the implicit
function theorem; we point out that the real transversal conditions yield then a real solution.
\end{enumerate}

If $\bx_i=(x_{i1},x_{i2},x_{i3})\in\R^3$ then $\bp_i=(p_{i1},p_{i2},p_{i3})\in\K^3$,
where $p_{ij}=(\xi_{ij}+O(t^{>0}))t^{-x_{ij}}$, $\xi_{ij}\ne0$ for all $1\le i\le N$, $j=1,2,3$. We denote
$\Ini(\bp_i)=\xi_i:=(\xi_{i1},\xi_{i2},\xi_{i3})\in(\C^*)^3$ and
$\Ini(\overline\bp)=\overline\xi:=(\xi_1,...,\xi_N)\in((\C^*)^3)^N$.

Introduce also the following auxiliary notation. If the circuit $C_S$ in the dual subdivision of $S$ is of type {\bf A}, we fix an
affine automorphism $\Phi_S:\Z^3\to\Z^3$ taking $C_S$ to a canonical pentatope $\Pi_{p,q}$
(see Section \ref{secA}). The
discriminantal equation of a polynomial $\sum_{\omega\in \Pi_{p,q}}a_\omega Z^\omega$ can be written in the form
\begin{equation}
(-1)^{1+p+q}a_{000}^{p+q}a_{100}^{-1}a_{010}^{-p}a_{001}^{-q}a_{1pq}=1\ .\label{emt5}
\end{equation}
Denote the exponent of a coefficient $a_\omega$ in this equation by $d(\omega)$.

\subsection{Singular points of tropical surfaces}\label{sing-trop}
By \cite[Theorem 2]{MMS}, the position of a singular point $\by\in S$ is defined uniquely
whenever the circuit $C_S$ is of type {\bf A}, {\bf B}, or {\bf D}. For circuit types
{\bf C} and {\bf E} there may be
several possible positions for $\by$. We will describe these possibilities via the geometry of
$\Graph(\nu_S)$. Namely, to determine the position of $\by$, it is enough to
determine the translation of $S$ which moves $\by$ to the origin. In turn, translations of
$S$ are in one-to-one correspondence with changes $\nu_S\mapsto\nu_S+\Lambda$, where
$\Lambda$ is any affine linear function. To move the singularity to the origin, we
use \cite[Lemma 10]{MMS}..

Without loss of generality, we assume that (cf. \cite[Theorem 2]{MMS})
$$C_S=\begin{cases}\{(1,0,0),(2,1,0),(0,2,0),(1,1,0)\},\quad &\text{if of type
{\bf C}},\\ \{(0,0,0),(0,0,1),(0,0,2)\},
\quad &\text{if of type {\bf E}}.\end{cases}$$

\begin{lemma}\label{lmt4}
 Let $C_S$ be of type {\bf C}, $\Lambda':\Delta\to\R$ the unique affine linear function, depending only on
 $x$ and $y$, which coincides with $\nu_S$ along $\conv(C_S)$. 
Set $\nu'=\nu_S-\Lambda'$ and
 introduce the following convex piecewise linear function
 on the projection $\pr_z(\Delta)$ of $\Delta$ to the $z$-axis: Set
$$\begin{cases}-c'_{\bm}=\min\{\nu'(\omega)\ :\ \omega\in\Delta\cap\Z^3,\ \pr_z(\omega)=\bm\},
\quad & \bm\in\pr_z(\Delta)\cap\Z\setminus\{0\},\\
-c'_0\gg\max\{-c'_{\bm},\ \bm\ne0\}, & \end{cases}$$
and then define a function $\nu_z:\pr_z(\Delta)\to\R$,
whose graph is the lower convex hull of $$\conv\left\{(\bm,-c'_{\bm})\ :\ \bm\in\pr_z(\Delta)\cap\Z
\right\}.$$ Then the possible singular points $\by\in S$ are in one-to-one correspondence with linear functions $\Lambda'':\Delta_z\to\R$ that vanish at the origin, are strictly less than
$\nu_z$, and whose graph is parallel to an edge
of $\Graph(\nu_z)$ which projects to one of the following segments:
\begin{equation}[-3,-1],\quad [-3,1],\quad [-3,3],\quad [-1,1],\quad [-1,3],\quad [1,3]\ .
\label{emt13}\end{equation}
\end{lemma}

{\bf Proof.}
The statement follows from \cite[Theorem 2 and Section 4.3]{MMS}:
According to the type of the weight class (see \cite[Lemma 10]{MMS}), we have to pick two points $\omega^{_1}$ and $\omega^{_2}$ in $\Delta\cap \Z^3$ whose coefficients $c_{\omega^{_1}}$ and $c_{\omega^{_2}}$ become equal and maximal among the $\omega\in \Delta\cap \Z^3 \setminus C_S$ after subtracting $\Lambda'$ and $\Lambda''$. Lemma 18 in \cite{MMS} yields the restriction that these points have to be picked with lattice distance one or three to the circuit $C_S$.
\proofend

\begin{lemma}\label{lmt5}
 Let $C_S$ be of type {\bf E}, $\Lambda':\Delta\to\R$ the unique affine linear function, depending only on
 $z$, which coincides with $\nu_S$ along $\conv(C_S)$. 
Set $\nu'=\nu_S-\Lambda'$ and
 introduce the following convex piecewise linear function
 on the projection $\pr_{x,y}(\Delta)$ of $\Delta$ to the $(x,y)$-plane: Set
$$\begin{cases}-c'_{\bm}=\min\{\nu'(\omega)\ :\ \omega\in\Delta\cap\Z^3,\ \pr_{x,y}(\omega)=\bm\},\quad
&\bm\in\pr_{x,y}(\Delta)\cap\Z^2\setminus\{0\},\\
-c'_0\gg\max\{-c'_{\bm},\ \bm\ne0\}, &\end{cases}$$
and then define a function $\nu_{x,y}:\pr_z(\Delta)\to\R$,
whose graph is the lower convex hull of $$\conv\left\{(\bm,-c'_{\bm})\ :\ \bm\in\pr_{x,y}(\Delta)\cap\Z^2\right\}.$$ Then the possible singular points $\by\in S$ are in one-to-one correspondence with
the linear functions $\Lambda'':\Delta_{x,y}\to\R$ that vanish at the origin, are strictly less than
$\nu_{x,y}$, and whose graph
\begin{enumerate}
\item[(i)] either is parallel to a triangular cell
of $\Graph(\nu_{x,y})$, whose
projection to the $(x,y)$-plane coincides up to
a $\Z$-linear transformation with one of the triangles:
\begin{equation}\begin{cases}&\conv\{(0,1),(1,0),(-1,-1)\},\quad \conv\{(0,1),(2,1),(-1,-1)\},\\
&\conv\{(0,1),(3,1),(-1,-1)\},\quad \conv\{(0,1),(3,1),(-3,-2)\},\\
&\conv\{(0,1),(4,1),(-2,-1)\},\quad\conv\{(-1,0),(0,1),(i,1)\},\ i\ge1\ .\end{cases}\label{emt6}\end{equation}
\item[(ii)] or is parallel to an edge $\widetilde E_1$ of
$\Graph(\nu_{x,y})$ and to a chord $\widetilde E_2$ joining
two vertices of $\Graph(\nu_{x,y})$ so that
\begin{itemize}
 \item projections of $\widetilde E_1,\widetilde E_2$ to the $(x,y)$-plane coincide up to
 a $\Z$-linear transformation with the pair
$$E_1=[(-1,0),(1,0)],\quad E_2=[(i,1),(j,-1)],\ i,j\in\Z\ ,$$
\item and the following condition holds:
\begin{eqnarray}&0<(\nu_{x,y}-\Lambda'')\big|_{E_1}<(\nu_{x,y}-\Lambda'')(\bm)
 \quad\text{for all}\ \bm\in\pr_{x,y}(\Delta)\cap\Z^2\setminus(E_1),\nonumber\\
 &(\nu_{x,y}-\Lambda'')\big|_{E_2}<(\nu_{x,y}-\Lambda'')(\bm)\quad\text{for all}\ \bm\in
 \pr_{x,y}(\Delta)\cap\Z^2\setminus(\Span(E_1)\cup E_2)\ .\nonumber
\end{eqnarray}
\end{itemize}
\end{enumerate}
\end{lemma}

{\bf Proof.}
The statement follows from \cite[Theorem 2, Propositions 21 and 23, and Section 4.6]{MMS}:
According to the type of the weight class (see \cite[Lemma 10]{MMS}), we have to pick either three points $\omega^{_1}$, $\omega^{_2}$ and $\omega^{_3}$ in $\Delta\cap \Z^3$ whose coefficients $c_{\omega^{_1}}$, $c_{\omega^{_2}}$ and $c_{\omega^{_3}}$ become equal and maximal among the $\omega\in \Delta\cap \Z^3 \setminus C_S$ after subtracting $\Lambda'$ and $\Lambda''$, or two pairs of points.

Let us first discuss the case of three points. Proposition 21 and Figure 17 in \cite{MMS} classify the possibilities up to $\Z$-linear transformation for the projections $\pr_{x,y}$ of the points $\omega^{_1}$, $\omega^{_2}$ and $\omega^{_3}$ under the assumption that there is no plane through $C_S$ such that they lie on the same side of this plane. This yields the first 5 possibilities of (\ref{emt6}).
The last case of (\ref{emt6}) is obtained if the points $\omega^{_1}$, $\omega^{_2}$ and $\omega^{_3}$ lie on the same side of a plane through $C_S$ following \cite[Proposition 23]{MMS}. Notice that the cases specified in \cite[Proposition 23(b,c)]{MMS} are not relevant.
Indeed, otherwise the function $\nu_{x,y}$ must be linear along a segment
containing at least three integral points. However this would yield that
the values of $\nu'$ at some three points $\bw_i,\bw_j,\bw_l$ outside $C_S$
are dependent with integral coefficients
which is impossible due to the
general choice of the values of $\nu'$ at these points (this generality results from
formulas (\ref{eA2}) and the generic choice of the parameters $M_i$ in (\ref{eCH})).

The case of a weight class with two pairs of points follows from \cite[Section 4.6]{MMS}.

Notice that every choice of $\Lambda''$ as described in the statement indeed yields a shift of $S$ with a tropically singular point at $0$: the vertices of the
triangles or pair of edges as specified above must satisfy certain arithmetic conditions (see \cite[Propositions 21 and 23, and Section 4.6]{MMS}).
We claim that these conditions are always satisfied. Indeed, these arithmetic restrictions geometrically mean that the convex hull
of the union of $C_S$ with the above points $\omega$ does not contain extra integral points.
However, if there were such an integral point, it would correspond to a vertex
of $\Graph(\nu_S)$, and this would break either the condition that $\nu_{x,y}$ is linear over the spoken triangle or edges, or
that $\Lambda''$ is strictly less than $\nu_{x,y}$.
\proofend

\subsection{Enhanced singular tropical surfaces}\label{enhanced}
Let us be given a point configuration $\overline\bx\in(\R^3)^N$ defined by (\ref{eCH}), a generic point configuration
$\overline\bp\in((\K^*)^3)^N$ such that $\Val(\overline\bp)=\overline\bx$, a tropical
surface $S\in\Sing^\tr(\Delta,\overline\bx)$ and its defining tropical
polynomial
\begin{equation}F_S(X)=\max_{\omega\in\Delta\cap\Z^3}(c_\omega+\langle X,\omega\rangle)\ .\label{etrpol}\end{equation}
Denote by $\nu_S:\Delta\to\R$ the convex piecewise linear function
Legendre dual to $F_S$, by $\Sigma_S$ the subdivision of $\Delta$ dual to
$S$, by $C_S$ the circuit, and by $P_S$ the corresponding
lattice path (formed by the edges dual to the $2$-faces of $S$ containing the points of $\overline\bx$).
Observe that $\nu_S(\omega)=-c_\omega$ for all points $\omega\in\Delta\cap\Z^3$.

\begin{lemma}\label{lmt2}
Any surface ${\mathcal S}\in\Sing(\Delta)$ that tropicalizes to $S$ is defined by a polynomial
\begin{equation}
\varphi_{\mathcal S}(\bz)=\sum_{\omega\in\Delta\cap\Z^3}(\alpha_\omega+O(t^{>0}))t^{\nu_S(\omega)}\bz^\omega
\in\K[\bz]\ ,
\label{emt1}\end{equation}
where $\bz^\omega=z_1^{\omega_1}z_2^{\omega_2}z_3^{\omega_3}$, and $O(t^{>0})$ accumulates the terms
containing $t$ to a positive power, and $\alpha_\omega\in\C^*$ for all $\omega
\in\Delta\cap\Z^3$. Furthermore, the polynomial
$$\Ini^{C_S}(\varphi_{\mathcal S})(Z):=\sum_{\omega\in C_S}\alpha_\omega Z^\omega
\in\C[Z]$$ has a singularity in $(\C^*)^3$.
\end{lemma}

{\bf Proof.}
We have to explain only the last claim. Viewing the surface ${\mathcal S}$ as an analytic equisingular family of
singular complex surfaces (cf. \cite[Section 2.3]{Sh}), we obtain an induced family of singular points with the limit
belonging to the big torus of $\Tor_\C(\delta)$ for some cell $\delta$ of the subdivision $\Sigma_S$,
that is, the cell dual to the face of $S$ containing the tropicalization of the singular point
(i.e., the tropical singular point). It is easy to see that,
for any cell $\delta\ne\conv(C_S)$ of $\Sigma_S$, any (nonzero) polynomial
$\sum_{\omega\in\delta\cap\Z^3}\beta_\omega Z^\omega$ has no singularities in $(\C^*)^3$.
Hence $\Ini^{C_S}(\varphi_{\mathcal S})$ must have singularity in $(\C^*)^3$.
\proofend

\begin{lemma}\label{lmt3}
If a polynomial $\varphi(\bz)$ of the form (\ref{emt1}) defines a surface in $(\K^*)^3$ passing through the configuration
$\overline\bp$, and if the polynomial $\Ini^{C_S}(\varphi)$ has a singularity in $(\C^*)^3$, then
the point $\overline\alpha:=(\alpha_\omega)_{\omega\in\Delta\cap \Z^3}\in\C
P^{N+1}$ belongs to a finite set
denoted by $A(S,\overline\bp)$. Furthermore,
\begin{enumerate}
 \item[(i)] If $C_S$ is of type {\bf A}, then
 \begin{itemize}
  \item for $P_S=\Gamma_k$, we have $|A(S,\overline\bp)|=|d(\Phi_S(\bw_k))|$;
  \item for $P_S=\Gamma_{k,k+1}$ and $\bw_{k+1}=\max C_S$, we have
  $|A(S,\overline\bp)|=|d(\Phi_S(\bw_{k+1}))|$;
  \item for $P_S=\Gamma_{k,k+1}$ and $\bw_{k+2}=\max C_S$, we have
  $$|A(S,\overline\bp)|=|d(\Phi_S(\bw_{k+2}))+d(\Phi_S(\bw_{k+1}))|\ .$$
\end{itemize}
 \item[(ii)] If $C_S$ is of type {\bf B}, then $|A(S,\overline\bp)|=\Vol_\Z(\conv(C_S))$
 when the tetrahedron $\conv(C_S)$ cannot be taken to
 $\conv\{(0,0,0),(1,0,0),(0,1,0),(3,7,20)\}$ by an automorphism of $\Z^3$ (cf. \cite[Theorem 2]{MMS}), and
 $|A(S,\overline\bp)|=\frac{1}{5}\Vol_\Z(\conv(C_S))=4$ when the tetrahedron $\conv(C_S)$ can be transformed to
 $\conv\{(0,0,0),(1,0,0),(0,1,0),(3,7,20)\}$ by an automorphism of $\Z^3$.
 \item[(iii)] If $C_S$ is of type {\bf C}, then $|A(S,\overline\bp)|=3$ (the lattice area of $\conv(C_S)$).
 \item[(iv)] If $C_S$ is of type {\bf D}, then $|A(S,\overline\bp)|=1$.
 \item[(v)] If $C_S$ is of type {\bf E}, then $|A(S,\overline\bp)|=1$ or $2$ according as $\conv(C_S)$ is an edge of the lattice path
 or not.
\end{enumerate}
\end{lemma}

{\bf Proof.}
We start by investigating the effect of the conditions imposed by the marked points $\bp_i$.
Tropically, the marked point $\bx_i$, $1\le i\le N$ lies on a $2$-face $F_i$ of $S$ dual to
an edge $E=[\omega^{_0},\omega^{_1}]\subset P_S$. In particular,
$$b:=c_{\omega^{_0}}+\langle\bx_i,\omega^{_0}\rangle=c_{\omega^{_1}}+\langle\bx_i,\omega^{_1}\rangle>
c_{\omega}+\langle\bx_i,\omega\rangle\quad\text{for all}\
\omega\in\Delta\setminus E\ ,$$ and then the condition imposed by the marked point $\bp_i$ is
$$0=\varphi(\bp_i)=t^{-b}\left(\Ini^E(\varphi)(\xi_i)+O(t^{>0})\right),
\quad\Ini^E(\varphi)(Z)=\sum_{\omega\in E}\alpha_\omega Z^\omega\ .$$
The lattice length $|E|:=|E\cap\Z^3|-1$ of $E$ is either $1$ or $2$. If
$|E|=1$, we obtain
\begin{equation}\alpha_{\omega^{_1}}=-\alpha_{\omega^{_0}}\xi_i^{\omega^{_0}-\omega^{_1}}\ .\label{emt2}\end{equation}
If $|E|=2$, then $\Ini^E(\varphi)(Z)$ has a singularity in $(\C^*)^3$; hence it is a monomial multiplied by
the square of a binomial,
which then implies
\begin{equation}\alpha_{\omega^{_1}}=\alpha_{\omega^{_0}}\xi_i^{\omega^{_0}-\omega^{_1}},\quad
\alpha_\omega=-2\alpha_{\omega^{_0}}
\xi_i^{(\omega^{_0}-\omega^{_1})/2},
\ \omega=\frac{\omega^{_0}+\omega^{_1}}{2}\\ .\label{emt3}\end{equation}

It follows, in particular, that $\overline\alpha$ is uniquely defined if $C_S$ is of type {\bf E} and
$\conv(C_S)$ is an edge of the lattice path $P_S$. If $C_S$ is of type {\bf E} and $\conv(C_S)\not\subset P_S$, then
we uniquely determine $\alpha_{\omega^{_0}}$ and $\alpha_{\omega^{_1}}$ for the end points $\omega^{_0},\omega^{_1}$ of $C_S$,
and by Lemma \ref{lmt2} obtain two values $\alpha_\omega=\pm2\sqrt{\alpha_{\omega^{_0}}\alpha_{\omega^{_1}}}$ for
the midpoint
$\omega$ of $C_S$, and hence two singular points of $\Ini^{C_S}(\varphi)$. Thus statement (v) is proved.

Now consider other types of circuits.

Suppose that $P_S=\Gamma_{k,k+1}$, $1\le k\le N$. As shown in Section \ref{sec-circuit},
$C_S$ must be of type {\bf A}. Equations (\ref{emt2}) yield $\overline\alpha\in\PP^{N+1}$ in the form
$$(\alpha_{\bw_0},...,\alpha_{\bw_k},\lambda\alpha'_{\bw_{k+1}},...,\lambda\alpha'_{\bw_{N+1}})\ ,$$
where
$(\alpha_{\bw_0},...,\alpha_{\bw_k},\alpha'_{\bw_{k+1}},...,\alpha'_{\bw_{N+1}})$ is a uniquely
defined generic point
of $\PP^{N+1}$, and $\lambda\ne0$ is an unknown parameter, which one can compute
from the discriminantal equation (\ref{emt5}) of the pentatope
$\Phi_S(\conv(C_S))$, obtaining
$|d(\Phi_S(\bw_k))|$ many solutions if $\bw_{k+1}=\max C_S$ and $|d(\Phi_S(\bw_{k+2}))+d(\Phi_S(\bw_{k+1}))|$ many solutions if $\bw_{k+2}=\max
C_S$.

Suppose that $P_S=\Gamma_k$, $0\le k\le N+1$. Then equations
(\ref{emt2}) determine the (nonzero) values
$\alpha_\omega$, $\omega\ne\bw_k$, up to proportionality.

If $C_S$ is of type {\bf A}, we obtain
$|d(\Phi_S(\bw_k))|$ values for the coefficient $\alpha_{\bw_k}$ from the discriminantal equation (\ref{emt5})
of the pentatope $\Phi_S(\conv(C_S))$. Thus (i) is proved.

If $C_S$ is of type {\bf B}, then $\bw_k$ is the interior point of the tetrahedron
$\conv(C_S)$. After a suitable transformation of the lattice $\Z^3$ and a coordinate change, we
obtain the equivalent question (see \cite[Theorem 2(a.2)]{MMS}):
How many values of $a\in\C^*$ are there such that the polynomial
$$\psi(x,y,z)=1+x+y+x^iy^jz^l+ax^{i'}y^{j'}z^{l'}$$ has a singularity in $(\C^*)^3$, where
$$(i,j,l)=(3,3,4),(2,2,5),(2,4,7),(2,6,11),(2,7,13),(2,9,17),(2,13,19),\ \text{or}\ (3,7,20)\ ,$$
and $(i',j',l')$ is the unique interior integral point of the tetrahedron $$\conv\{(0,0,0),(1,0,0),(0,1,0),(i,j,l)\}?$$
The system of equations $\psi=\psi_x=\psi_y=\psi_z=0$ reduces to
\begin{equation}x=\lambda,\quad y=\mu,\quad z^l=\nu,\quad a=\rho z^{-l'}\label{emt11}\end{equation} with
some nonzero constants $\lambda,\mu,\nu,\rho$. In all cases except for
$(i,j,l)=(3,7,20)$, we have $\gcd(l',l)=1$, and hence
$l=\Vol_\Z(\conv(C_S))$ solutions for $a$. In the remaining case
$(i,j,k)=(3,7,20)$, $(i',j',l')=(1,2,5)$, and we obtain $4=
\Vol_\Z(\conv(C_S))/5$
values for $a$. The proof of (ii) is completed.

If $C_S$ is of type {\bf C}, then $\bw_k$ is the interior point of the triangle
$\conv(C_S)$. For given $\alpha_\omega\ne0$ at the vertices $\omega$ of $\conv(C_S)$,
there are exactly $\Vol_\Z(\conv(C_S))=3$ values $\alpha_{\bw_k}$, corresponding to singular
polynomials $\Ini^C(\varphi_{\mathcal S})$ (cf. \cite[Lemma 3.5]{Sh}). This yields (iii).

If $C_S$ is of type {\bf D}, then $\bw_k$ is a vertex of the parallelogram
$\conv(C_S)$. If $\bw_i,\bw_j,\bw_l$ are the other vertices of $\conv(C_S)$, and $\bw_j$
is opposite to $\bw_k$, then the fact that
$\Ini^{C_S}(\varphi_{\mathcal S})$ has a singularity in $(\C^*)^3$ yields
$\alpha_{\bw_k}=\alpha_{\bw_i}\alpha_{\bw_l}\alpha_{\bw_j}^{-1}$, which defines $\alpha_{\bw_k}$ uniquely. Thus statement (iv) is proved. \proofend

\begin{remark}\label{rmt1}
 Observe that, in the case of the lattice path $\Gamma_{k,k+1}$ and a circuit of type
 {\bf A} containing the points $\bw_{k+1},\bw_{k+2}$, one may obtain an empty set
 $A(S,\overline\bp)$.
\end{remark}

We call the points $\overline\alpha\in A(S,\overline\bp)$ \emph{enhancements} of $S$, and the pairs
$(S,\overline\alpha)$ \emph{enhanced singular tropical surfaces}.

\subsection{Patchworking construction}
\label{sec-pw}
We will now see how given enhancements $\overline\alpha$ can be lifted
to equations of algebraic surfaces in $\Sing(\Delta,\overline\bp,S)$
using patchworking techniques. The idea is to look for a solution in the form
\begin{equation}\varphi_{\mathcal S}(\bz)=\sum_{\omega\in\Delta\cap\Z^3}a_\omega t^{\nu_S(\omega)}\bz^\omega\in\K[\bz]
\label{erev1}\end{equation} (cf. formula (\ref{emt1})), where  $a_{\omega^{_0}}\equiv1$ for some vertex $\omega^{_0}$ of the subdivision
$\Sigma_S$, and the remaining coefficients $a_\omega=\alpha_\omega+O(t{^>0})$ are obtained from the conditions to
pass through $\overline\bp$ and to have a singular point $\bq$ with $\Ini(\bq)=z$, a singular point
of  $\Ini^{C_S}(\varphi_{\mathcal S})(Z)=\sum_{\omega\in C_S}\alpha_\omega Z^\omega$ in $(\C^*)^3$.
We then show that, in the case of circuits of type
{\bf A}, {\bf B}, {\bf C}, and {\bf D},
at $t=0$, these conditions turn into a system of equations with a non-degenerate linearization, and then apply the implicit function theorem.
In case of circuits of type {\bf E}, sometimes one has to use additional terms in the Puiseux series representing
the coordinates of the points of $\overline\bp$ in order to get a non-degenerate system of equations.

In Lemmas \ref{lmt6}-\ref{lmt7}, we settle the case of circuits of type {\bf A}, {\bf B}, {\bf C}, and {\bf D},
and in Lemma \ref{lmt6} the most difficult case,
namely circuits of type {\bf E}.

\begin{lemma}\label{lmt6}
Let $S\in\Sing^\tr(\Delta,\overline\bx)$, and let the circuit $C_S$ in the subdivision dual to $S$ be of type {\bf A} or {\bf B}
(see lemma \ref{lrev1}). Then
\begin{enumerate}
\item[(i)] if $C_S$ is not $\Z$-affine equivalent to
$\{(0,0,0),(1,0,0),(0,1,0),(3,7,20)\}$, then, for any
point $\overline\alpha\in A(S,\overline\bp)$, (see Lemma \ref{lmt3}) there exists a unique
algebraic surface ${\mathcal S}\in\Sing(\Delta,\overline\bp,S)$;
\item[(ii)] if $C_S$ is $\Z$-affine equivalent to
$\{(0,0,0),(1,0,0),(0,1,0),(3,7,20)\}$, then, for any
point $\overline\alpha\in A(S,\overline\bp)$, there exist $5$ algebraic surfaces
${\mathcal S}\in\Sing(\Delta,\overline\bp,S)$ matching the enhancement $\overline\alpha$.
\end{enumerate}
\end{lemma}

{\bf Proof.} The required statement can again be viewed as a patchworking theorem, and it follows from
a suitable version of the implicit function theorem. Namely, we look for polynomials given
by (\ref{erev1})
(cf. formula (\ref{emt1})), where  $a_{\omega^{_0}}\equiv1$ for some vertex $\omega^{_0}$ of the subdivision
$\Sigma_S$, and the remaining coefficients $a_\omega=\alpha_\omega+O(t{^>0})$ are obtained from the conditions to
pass through $\overline\bp$ and to have a singular point $\bq$ with $\Ini(\bq)=z$, a singular point
of  $\Ini^{C_S}(\varphi_{\mathcal S})(Z)=\sum_{\omega\in C_S}\alpha_\omega Z^\omega$ in $(\C^*)^3$. At $t=0$ these conditions turn into the system
of equations (\ref{emt2}) in the coefficients $a_\omega$, $\omega\neq \omega^{_0}$, and the discriminantal equation for the circuit $C_S$.

In the case (i), if the lattice path is $\Gamma_k$ for some $k$, the Jacobian of the above system at $t=0$ is a (suitably arranged) lower triangular matrix with
the nonzero
entries from (\ref{emt2}) and the discriminantal equation for the circuit. If the lattice path is $\Gamma_{k,k+1}$, we obtain a matrix whose column corresponding to $a_{\bw_{k+1}}$ has only one nonzero entry, and is (suitably arranged) lower triangular after erasing this column and the corresponding row.

In the case (ii), without loss of generality suppose that
$$C_S=\{(0,0,0),(1,0,0),(0,1,0),(3,7,20),(1,2,5)\}\ .$$ Any singular complex polynomial $\Ini^{C_S}(\varphi_{\mathcal S})$ supported at
$C_S$ has
$5$ singular points in $(\C^*)^3$, obtained from each other by
the $\Z_5$-action $z\mapsto z\eps$, $\eps^5=1$. Each singular point $z\in(\C^*)^3$ of $\Ini^{C_S}(\varphi_{\mathcal S})$ is an ordinary node,
in particular,
\begin{equation}\det(\Hessian(\Ini^{C_S}(\varphi_{\mathcal S}))(z))\ne0\ .\label{emt9}\end{equation}
Then we consider the system of equations in the coefficients $a_\omega$, $\omega\ne\omega^{_0}$, of the sought polynomial
$\varphi_{ \mathcal S}$ and the coordinates $z_i$, $i=1,2,3$, of the singular point, where $z_i=z_{i0}+O(t^{>0})$. The equations induced by the conditions
${\mathcal S}\supset\overline\bp$ and
\begin{equation}\varphi_{\mathcal S}(z_1,z_2,z_3)=\frac{\partial\varphi_{\mathcal S}}{\partial z_i}(z_1,z_2,z_3)=0,\ i=1,2,3\ ,
\label{emt10}\end{equation}
and this system has a unique solution, since its Jacobian at $t=0$ does not vanish:
\begin{itemize}\item the block coming
from the conditions $\varphi_{\mathcal S}(\bp_i)=0$, $i=1,...,N$, is the Jacobian of the nondegenerate linear system
(\ref{emt2}), \item for the block coming from the system  (\ref{emt10}), the nondegeneracy follows from (\ref{emt9}).\proofend
\end{itemize}

\begin{lemma}\label{lmt8}
Let $S\in\Sing^\tr(\Delta,\overline\bx)$, $C_S$ be of type {\bf C}. Let us be given an enhancement
$\overline\alpha\in A(S,\overline\bx)$ and a tropical singular point $\by\in S$ associated with
a segment
$\sigma=[m,n]$ as specified in Lemma \ref{lmt4}, formula (\ref{emt13}). Then there are
$(n-m)$ algebraic surfaces ${\mathcal S}\in\Sing(\Delta,\overline\bp,S)$,
 matching the given data $\overline\alpha$ and $\by$.
\end{lemma}

{\bf Proof.} Without loss of generality we can suppose that the lattice path $P_S=\Gamma_k$, the left out point $\bw_k$ is $(1,1,0)$, the circuit is
$C_S=\{(1,0,0),(2,1,0),(0,2,0),(1,1,0)\}$, 
the tropical singular point
is $\by=(0,0,0)$, and the sought polynomial
takes form (cf. \cite[Theorem 2(b.1)]{MMS})
\begin{eqnarray}\varphi_S(\bz)&=&\sum_{(i,j,0)\in C_S}a_{ij0}z_1^iz_2^j+\sum_{(i,j,0)\in\Delta\setminus C_S}O(t^{>0})\cdot
z_1^iz_2^j\nonumber\\
& &+t^s\left(a_{i_1j_1m}z_1^{i_1}z_2^{j_1}z_3^m+a_{i_2j_2n}z_1^{i_2}z_2^{j_2}z_3^n\right)
+O(t^{>s}),\nonumber\end{eqnarray} where $s>0$, and
$$a_{100}\equiv1,\ a_{210}=\alpha_{210}+O(t^{>0}),\ a_{020}=\alpha_{020}+O(t^{>0}),\ a_{110}
=\alpha_{110}+O(t^{>0}),$$
$$a_{i_1j_1m}=\alpha_{i_1j_1m}+O(t^{>0}),\quad a_{i_2j_2n}=\alpha_{i_2j_2n}+O(t^{>0})\ .$$
The equations
\begin{align*}(\varphi_S)_{t=0}(z_{10},z_{20},z_{30})=&\left(\frac{\partial\varphi_S}{\partial z_1}\right)_{t=0}(z_{10},z_{20},z_{30})\\=&
\left(\frac{\partial\varphi_S}{\partial z_2}\right)_{t=0}(z_{10},z_{20},z_{30})=\left(t^{-s}\frac{\partial\varphi_S}{\partial z_3}\right)_{t=0}
(z_{10},z_{20},z_{30})=0\end{align*} for $(z_{10},z_{20},z_{30})=\Ini(\bq)$, $\bq$ being a singular point of the sought surface ${\mathcal S}\in\Sing(\Delta,\overline\bp,S)$,
give a unique solution $(z_{10},z_{20})$ for the singularity of $\Ini^{C_S}(\varphi_S)$ in $(\C^*)^2$, and the last equation,
\begin{equation}m\alpha_{i_1j_1m}z_{10}^{i_1}z_{20}^{j_1}z_{30}^{m-1}+n\alpha_{i_2j_2n}
z_{10}^{i_2}z_{20}^{j_2}
z_{30}^{n-1}=0\ ,
\label{emt14}\end{equation}
yields $(n-m)$ nonzero solutions for $z_{30}$.

We claim that each solution induces a unique
surface ${\mathcal S}\in\Sing(\Delta,\overline\bp,S)$ matching the requirements of the lemma.
Indeed, the implicit function theorem applies: the system
$$\varphi_S(\bp_i)=0,\quad i=1,...,N\ ,$$
linearizes into the nondegenerate linear system
(\ref{emt2}) with respect to the variables $a_\omega$, $\omega\in\Delta\cap\Z^3\setminus\{(1,0,0),(1,1,0)\}$,
and the Jacobian evaluated at $t=0$ for the system
$$\varphi_S(\bq)=\frac{\partial\varphi_S}{\partial z_1}(\bq)=\frac{\partial\varphi_S}{\partial z_2}(\bq)=t^{-s}
\frac{\partial\varphi_S}{\partial z_3}(\bq)=0$$ with respect to $a_{110}$ and the coordinates of $\bq$
takes form of a lower block-triangular matrix
$$\left(\begin{matrix} z_{10}z_{20} & 0 & 0\\
* & \Hessian(\Ini^{C_S}(\varphi))(z_{10},z_{20}) & 0\\ * & * & Q_{zz}(z_{10},z_{20},z_{30}) \end{matrix}\right)$$
where $Q=\alpha_{i_1j_1m}z_1^{i_1}z_2^{j_1}z_3^m+\alpha_{i_2j_2n}z_1^{i_2}z_2^{j_2}z_3^n$.
The nondegeneracy of this matrix (coming particularly from the fact that $(z_{10},z_{20})$ is
an ordinary node of $\Ini^{C_S}(\varphi)$ and that the nonzero roots $z_{30}$ of (\ref{emt14}) are simple) completes the proof.
\proofend

\begin{lemma}\label{lmt7}
Let $S\in\Sing^\tr(\Delta,\overline\bx)$, $C_S$ be of type {\bf D}. Then there are two surfaces
${\mathcal S}\in\Sing(\Delta,\overline\bp,S)$.
\end{lemma}

{\bf Proof.}
Without loss of generality we can suppose that the lattice path $P_S=\Gamma_k$, the left out point $\bw_k$ is the origin, the circuit is
$C_S=\{(0,0,0),(1,0,0),(0,1,0),(1,1,0)\}$, 
the unique tropical singular point
is $\by=(0,0,0)$, and the sought polynomial
takes the form (cf. \cite[Theorem 2(b.2)]{MMS})
\begin{eqnarray}\varphi_S(\bz)&=&z_1z_2+a_{100}z_1+a_{010}z_2+a_{000}+
\sum_{(u,v,0)\in\Delta\setminus C_S}
O(t^{>0})\cdot z_1^uz_2^v\nonumber\\
& &+t^s(a_{ij1}z_1^i z_2^j
z_3+a_{m,n,-1}z_1^my^nz_3^{-1})+O(t^{>s})\ ,\nonumber\end{eqnarray}
where $s>0$, and
$$a_{100}=\alpha_{100}+O(t^{>0}),\ a_{010}=
\alpha_{010}+O(t^{>0}),\ a_{000}=\alpha_{000}+
O(t^{>0})\ ,$$ $$a_{ij1}=\alpha_{ij1}+O(t^{>0}),
\ a_{m,n,-1}=\alpha_{m,n,-1}+O(t^{>0})\ .$$ A possible singular point of a sought surface ${\mathcal S}\in\Sing(\Delta,\overline\bp,S)$
should be
$\bq=(z_{10}+O(t^{>0}),z_{20}+O(t^{>0}),z_{30}+O(t^{>0}))$, where $(z_{10},z_{20},z_{30})$ are found from the system
$\varphi\big|_{t=0}=\frac{\partial \varphi}{\partial
  z_1}\big|_{t=0}=\frac{\partial \varphi}{\partial
  z_2}\big|_{t=0}=(t^{-s}\frac{\partial\varphi}{\partial z_3})\big|_{t=0}=0$, which reduces to
\begin{equation}
\alpha_{000}=\alpha_{100}\alpha_{010},\quad z_{10}=-\alpha_{010},\quad z_{20}=-\alpha_{100},\quad
\alpha_{ij1}z_{10}^i
z_{20}^j-\alpha_{m,n,-1}z_{10}^mz_{20}^nz_{30}^{-2}=0\ .
\label{ereal1}\end{equation}
Thus, we get two solutions $(z_{10},z_{20},z_{30})$, and we claim that each of them
induces a unique algebraic surface ${\mathcal S}\in\Sing(\Delta,\overline\bp,S)$. Again we apply the implicit function theorem to the system
of equations
\begin{equation}\varphi(\bp_i)=0,\ i=1,...,N,\quad \varphi(\bq)=\frac{\partial\varphi}{\partial z_1}(\bq)=
\frac{\partial\varphi}{\partial z_2}(\bq)=t^{-s}\frac{\partial\varphi}{\partial z_3}(\bq)=0\label{emt12}
\end{equation}
in the coordinates of $\bq$ and the coefficients $a_\omega$, $\omega\in\Delta\cap\Z^3\setminus\{(1,1,0)\}$.
Similarly to the proof of Lemma \ref{lmt6}, the Jacobian, evaluated at $t=0$, splits into a block coming from the
nondegenerate linear system (\ref{emt2}) and a block coming from the last four equations in (\ref{emt12}):
$$\left(\begin{matrix} 1 & 0 & 0 & 0\\ 0 & \alpha_{100} & 0 & 0\\
0 & 0 & \alpha_{010} & 0\\ 0 & * & * &
2\alpha_{m,n,-1}z_{10}^mz_{20}^nz_{30}^{-3}\end{matrix}\right)$$ that is nondegenerate too.
\proofend

\begin{lemma}\label{lmt11}
Let $S\in\Sing^\tr(\Delta,\overline\bx)$ have a circuit $C_S$ of type {\bf E}, and
$C_S=\{(0,0,0),(0,0,1),(0,0,2)\}$. 

(1) Suppose that a singular point $\by\in S$ is associated with
a triangle $\delta\subset\R^2\hookrightarrow\R^3$ from the list
(\ref{emt6}), as specified in Lemma \ref{lmt5}(i). Then there exist precisely
$2\cdot\Vol_\Z(\delta)$ surfaces ${\mathcal S}\in\Sing(\Delta,\overline\bp,S)$ that have a singular point tropicalizing to $\by$.

(2) Suppose that a singular point $\by\in S$ is associated with
a pair of edges $E_1,E_2$ as specified in Lemma \ref{lmt5}(ii). Then there exist
precisely
$8$ surfaces ${\mathcal S}\in\Sing(\Delta,\overline\bp,S)$ that have a singular point tropicalizing to $\by$.
\end{lemma}

{\bf Proof.} In both cases the lattice path is $P_S=\Gamma_k$ for some $1\le k\le N$. Furthermore, we have the following options:
\begin{enumerate}\item[(i)] either the segment $\conv(C_S)$ is a part of the lattice path $\Gamma_k$, and its dual $2$-face of
$S$ contains a marked point $\bx_{k_0}=(\lambda,\mu,0)$, where we can suppose that
$\lambda,\mu$ are generic in the sense of (\ref{e6:6});
\item[(ii)] or $\conv(C_S)$ is not an edge of $\Gamma_k$.
\end{enumerate}

\smallskip

{\it Step 1.} Consider the possibility (i).
Then the enhancement $\overline\alpha$
is uniquely restored from formulas (\ref{emt2}) and (\ref{emt3}), when we set $\omega^{_0}=(0,0,2)$.
We have
\begin{equation}
  \varphi_S(\bz)=z_3^2-2a_{001}z_3+a_{000}+\sum_{\omega\in\Delta\cap\Z^3\setminus C_S}a_\omega t^{-c_\omega} \bz^\omega
  \label{emt17}
\end{equation}
with $-c_\omega>0$ for all $\omega\in\Delta\cap\Z^3\setminus C_S$ and
$\Val(a_\omega)=0$ for all $\omega$, and
we have
\begin{displaymath}
  \Delta\cap\Z^3=\{w_0,\ldots,w_{N+1}\}
\end{displaymath}
with
\begin{displaymath}
  C_S=\{w_{k-1}=(0,0,2),w_k=(0,0,1),w_{k+1}=(0,0,0)\}.
\end{displaymath}
We intend to solve the system of equations
$$\varphi_S(\bp_i)=0,\ i=1,...,N,\quad\varphi_S(\bq)=\frac{\partial\varphi_S}{\partial x}(\bq)=
\frac{\partial\varphi_S}{\partial y}(\bq)=\frac{\partial\varphi_S}{\partial z}(\bq)=0$$
with respect to the variables $a_{001}$, $a_{000}$, $a_\omega$, $\omega\in\Delta\cap\Z^3\setminus C_S$,
and the coordinates $z_1,z_2,z_3$ of the singular point $\bq$ with the
aid of the implicit function theorem.

Recall that, in the framework
of Lemma \ref{lmt5}, $\by=\Val(\bq)$ is the origin, i.~e.~$\Val(z_i)=0$, and
let $\Ini(\bq)=(z_{10},z_{20},z_{30})$.
Indeed,
$z_{30}=1$, which follows from the equation $\varphi_S(\bq)=0$ and
(\ref{emt17}).

(1) In the first case, let $\delta=\conv\{(i_1,j_1),(i_2,j_2),(i_3,j_3)\}$. There exist uniquely
defined $l_1,l_2,l_3\in\Z$ such that, in the notation of Lemma \ref{lmt5},
$$(i_r,j_r,l_r)\in\Delta\quad\text{and}\quad\nu_{x,y}(i_r,j_r)=\nu'(i_r,j_r,l_r),\ r=1,2,3\ .$$
Then by formula (\ref{emt17}) we have
\begin{displaymath}
  -c_{i_rj_rl_r}=s<-c_\omega
\end{displaymath}
for $r=1,2,3$ and all other $\omega\not\in\{\omega_1=\omega_2=0\}$.
The equations
$$\left(t^{-s}\frac{\partial\varphi_S}{\partial z_1}\right)_{t=0,z_3=z_{30}}=
\left(t^{-s}\frac{\partial\varphi_S}{\partial z_2}\right)_{t=0,z_3=z_{30}}=0$$
yield that the coordinates $z_{10},z_{20}$ of $\Ini(\bq)=(z_{10},z_{20},z_{30})$ correspond to critical points
in $(\C^*)^2$ of
the polynomial
\begin{displaymath}
  Q(z_1,z_2)=\sum_{r=1}^3\alpha_{i_rj_rl_r}z_1^{i_r}z_2^{j_r},
\end{displaymath}
which gives us $\Vol_\Z(\delta)$
solutions $(z_{10},z_{20})$ as possible initial values for $\bq$ in total.

In order to apply the implicit function theorem we have to replace the
equation $\varphi_S(\bq)=0$ by two possible other equations, since it is
unsuitable itself being of degree two in $z_3$. We now first want to
derive these new equations. For that we consider the equation $\varphi_S(\bp_k)=0$,
\begin{equation}
  1-2a_{001}+a_{000}+\sum_{\omega\in\Delta\cap\Z^3\setminus
    C_S}a_\omega t^{-c_\omega+\lambda\omega_1+\mu\omega_2}=0,
  \label{e6}\nonumber
\end{equation}
together with $\frac{\partial\varphi_S}{\partial z_3}(\bq)=0$,
\begin{equation}
  \label{e6:1}\nonumber
  2z_3-2a_{001}+
  \sum_{\omega\in\Delta\cap\Z^3\setminus C_S}a_\omega
  t^{-c_\omega}\omega_3 z_1^{\omega_1}z_2^{\omega_2}z_3^{\omega_3-1}=0.
\end{equation}
The equations lead to
\begin{equation}
  \label{e6:2}\nonumber
  a_{001}=z_3+\frac{1}{2}\cdot\sum_{\omega\in\Delta\cap\Z^3\setminus C_S}a_\omega
  t^{-c_\omega}\omega_3 z_1^{\omega_1}z_2^{\omega_2}z_3^{\omega_3-1}
\end{equation}
and
\begin{equation}
  \label{e6:3}\nonumber
  a_{000}=-1+2z_3+\sum_{\omega\in\Delta\cap\Z^3\setminus C_S}a_\omega\cdot
  \big(t^{-c_\omega}\omega_3 z_1^{\omega_1}z_2^{\omega_2}z_3^{\omega_3-1}-t^{-c_\omega+\lambda\omega_1+\mu\omega_2}\big).
\end{equation}
Plugging these equations into (\ref{emt17}) and reorganizing the terms
we get
\begin{equation}
  \label{e6:4}\nonumber
  (z_3-1)^2=\sum_{\omega\in\Delta\cap\Z^3\setminus C_S}a_\omega t^{-c_\omega}\cdot
  \big(Z^\omega-t^{\lambda\omega_1+\mu\omega_2}
  +(z_3-1)\cdot \omega_3 z_1^{\omega_1}z_2^{\omega_2}z_3^{\omega_3-1}\big),
\end{equation}
and taking square roots we get two equations
\begin{equation}
  \label{e6:5}
  \psi_{\pm}=z_3-1\pm\sqrt{\sum_{\omega\in\Delta\cap\Z^3\setminus C_S}a_\omega t^{-c_\omega}\cdot
  \big(Z^\omega-t^{\lambda\omega_1+\mu\omega_2}
  +(z_3-1)\cdot \omega_3 z_1^{\omega_1}z_2^{\omega_2}z_3^{\omega_3-1}\big)}=0
\end{equation}
to replace $\varphi_S(\bq)=0$ with.

We now consider the polynomial map $\Psi$, that maps
\begin{displaymath}
  \zeta=(t,a_{w_1},\ldots,a_{w_{k-2}},z_1,z_2,z_3,a_{w_k},\ldots,a_{w_{N+1}})
\end{displaymath}
to
\begin{displaymath}
  \left(t^{-s_1}\varphi_S(\bp_1),\ldots,t^{-s_{k-1}}\varphi_S(\bp_{k-1}),
    t^{-s}\frac{\partial\varphi_S}{\partial z_1},t^{-s}\frac{\partial\varphi_S}{\partial z_2},
    \psi_{\pm},\frac{\partial\varphi_S}{\partial z_3},t^{-s_{k}}\varphi_S(\bp_{k}),\ldots,
    t^{-s_{N-1}}\varphi_S(\bp_{N-1})\right),
\end{displaymath}
where $s_i=\Val(\varphi_S(\bp_i))$. Note, that the initial values give a zero
\begin{displaymath}
  \zeta_0=(0,\alpha_{w_1},\ldots,\alpha_{w_{k-2}},z_{10},z_{20},z_{30},\alpha_{w_k},\ldots,\alpha_{w_{N+1}})
\end{displaymath}
of $\Psi$. We assume that the values $\lambda$ and $\mu$ are generic
in the sense that
\begin{equation}
  -c_{\omega'}+\lambda\omega_1'+\mu\omega_2'\not=
  -c_{\omega}+\lambda\omega_1+\mu\omega_2\not=-c_{\omega'}
  \label{e6:6}
\end{equation}
whenever $\omega\not=\omega'$. This ensures that the term under the square
root in $\psi_\pm$ is non-zero, if we evaluate it at
$(\alpha_{w_1},\ldots,\alpha_{w_{k-2}},z_{10},z_{20},z_{30},\alpha_{w_k},\ldots,\alpha_{w_n})$,
so that $\psi_\pm$ is analytic locally in $\zeta_0$. Moreover,
computing derivatives in (\ref{e6:5})  we get
\begin{displaymath}
  \frac{\partial\psi_\pm}{\partial z_3}(\zeta_0)=1\pm
  \frac{
    \sum_{\omega\in\Delta\cap\Z^3\setminus C_S}a_\omega t^{-c_\omega}\cdot
    (\omega_3z_1^{\omega_1}z_2^{\omega_2}z_3^{\omega_3-1}+(z_3-1)\omega_3(\omega_3-1)z_1^{\omega_1}
    z_2^{\omega_2}z_3^{\omega_3-2})
    }{
      2\sqrt{
        \sum_{\omega\in\Delta\cap\Z^3\setminus C_S}a_\omega t^{-c_\omega}\cdot
        (Z^\omega-t^{\lambda\omega_1+\mu\omega_2}
        +(z_3-1)\cdot \omega_3 z_1^{\omega_1}z_2^{\omega_2}z_3^{\omega_3-1})
      }
    }_{\big|\zeta=\zeta_0}=1
\end{displaymath}
and all other derivatives of $\psi_\pm$ vanish at $\zeta_0$, since due
to the genericity assumption on $\lambda$ and $\mu$ the valuation of
the denominator in the fraction is at most half the valuation of the
numerator. Similar computations for the other component functions of
$\Psi$ lead to
the following Jacobian of $\Psi$ with respect to all variables but $t$
evaluated at $\zeta_0$,
\begin{displaymath}
  \left(
  \begin{array}{cccc|cccccccc}
    \xi_1^{w_0}&\xi_1^{w_1}&&\\
    &\ddots&\ddots&\\
    &&\xi_{k-1}^{w_{k-2}}&\xi_{k-1}^{w_{k-1}}\\\hline
    *&\ldots&\ldots&*&\frac{\partial^2 t^{-s}\varphi_S}{\partial z_1^2}(\zeta_0)&\multicolumn{1}{c|}{\frac{\partial^2 t^{-s}\varphi_S}{\partial z_2\partial z_1}(\zeta_0)}&*&\ldots&\ldots&\ldots&\ldots&*\\
    *&\ldots&\ldots&*&\frac{\partial^2 t^{-s}\varphi_S}{\partial z_1\partial z_2}(\zeta_0)&\multicolumn{1}{c|}{\frac{\partial^2 t^{-s}\varphi_S}{\partial z_2^2}(\zeta_0)}&*&\ldots&\ldots&\ldots&\ldots&*
  \\\cline{5-12}
    &&&&&\multicolumn{1}{c|}{}&1&&&&&\\
    &&&&&\multicolumn{1}{c|}{}&2&-2&&&&\\
    &&&&&\multicolumn{1}{c|}{}&&-2&1&&&\\
    &&&& &\multicolumn{1}{c|}{} & &  &\xi_{k+1}^{w_{k+1}}&\xi_{k+1}^{w_{k+2}}&&\\
    &&&& &\multicolumn{1}{c|}{} & &  &&\ddots&\ddots&\\
    &&&& &\multicolumn{1}{c|}{} & &  &&&\xi_{N-1}^{w_{N-1}}&\xi_{N-1}^{w_N}
  \end{array}
  \right),
\end{displaymath}
where all missing entries are zero and the stars denote possibly
non-zero entries. Since the critical point $(z_{10},z_{20})$ of $Q$
is non-degenerate the Hessian in the middle block has a non-vanishing
determinant and thus the determinant of the Jacobian does not
vanish. Applying the implicit function theorem we get in each of the
two cases $\psi_-$ and $\psi_+$ a unique solution, and since we have
$\Vol_\Z(\delta)$ choices for $(z_{10},z_{20})$ we end up with
$2\cdot\Vol_\Z(\delta)$ algebraic surfaces ${\mathcal S}\in\Sing(\Delta,
\overline\bp,S)$ having a singular point $\bq$ with $\Trop(\bq)=\by$.

\smallskip

(2) The second case works along the same lines. With the notation of
Lemma \ref{lmt5} the relations
$$\begin{cases}& (-1,0,l_1),(1,0,l_2),(i,1,l_3),(j,-1,l_4)\in\Delta\cap\Z^3,\\
&\nu_{x,y}(-1,0)=\nu'(-1,0,l_1),\ \nu_{x,y}(1,0)=\nu'(1,0,l_2),\\
&\nu_{x,y}(i,1)=\nu'(i,1,l_3),\ \nu_{x,y}(j,-1)=\nu'(j,-1,l_4)\ \end{cases}$$
uniquely determine integers $l_1,l_2,l_3,l_4$ and valuations
$s_2>s_1>0$, such that
\begin{displaymath}
  s_1=-c_{-1,0,l_1}=-c_{1,0,l_2}<-c_\omega
\end{displaymath}
for all other $\omega\in\Delta\cap\Z^3$ of the form $\omega=(i,0,l)$
and such that
\begin{displaymath}
  s_2=-c_{i,1,l_3}=-c_{j,1,l_4}<-c_\omega
\end{displaymath}
for all remaining $\omega\not\in\{\omega_2=0\}$. Defining
$$Q_1(z_1,z_3)=\alpha_{-1,0,l_1}z_1^{-1}z_3^{l_1}+\alpha_{1,0,l_2}z_1z_3^{l_2},\quad
Q_2(z_1,z_2,z_3)=\alpha_{i,1,l_3}z_1^iz_2z_3^{l_3}+\alpha_{j,1,l_4}z_1^jz_2^{-1}z_3^{l_4}\ ,$$
the critical points of $Q_1$ and $Q_2$ respectively determine the
possible pairs $(z_{10},z_{20})$ for $\Ini(\bq)=(z_{10},z_{20},1)$ via
the equations $$\left(t^{-s_1}\frac{\partial\varphi_S}{\partial z_1}\right)_{t=0}(z_{10},z_{20},1)=\left(t^{-s_2}\frac{\partial\varphi_S}{\partial z_2}
\right)_{t=0}(z_{10},z_{20},1)=0.$$ They are thus the solutions of the system
\begin{equation}-\alpha_{-1,0,l_1}z_{10}^{-2}+\alpha_{1,0,l_2}=\alpha_{i,1,l_3}
z_{10}^i-\alpha_{j,-1,l_4}z_{10}^jz_{20}^{-2}=0\ ,
\label{e3}\end{equation}
from which we get $4$ solutions
$(z_{10},z_{20})\in(\C^*)^2$. Replacing the equations
$t^{-s}\frac{\partial \varphi_S}{\partial z_1}$ and
$t^{-s}\frac{\partial \varphi_S}{\partial z_2}$ in case (1) by
$t^{-s_1}\frac{\partial \varphi_S}{\partial z_1}$ and
$t^{-s_2}\frac{\partial \varphi_S}{\partial z_2}$, we can continue as
in case (1) and find $8$ surfaces ${\mathcal S}\in\Sing(\Delta,\overline\bp,S)$ having a singular point
$\bq$ tropicalizing to $\by$.

\smallskip

{\it Step 2.} In the situation (ii), the above argument appears to be rather simpler.
Note that the equations $\varphi_S(\bp_i)=0$, $i=1,...,N$, and the condition
$\omega^{_0}=(0,0,2)$ uniquely determine the values
$\alpha_\omega$ for all $\omega\ne\bw_k$. For $\omega=\bw_k$, we obtain two values
\begin{equation}
\alpha_{\bw_k}=\alpha_{001}=\pm\sqrt{\alpha_{000}}
\label{e-real2}\end{equation} (cf. Lemma \ref{lmt2}), and respectively
$z_{30}=-\alpha_{001}$. Independently of the
choice of $\alpha_{001}$, we obtain $\Vol_\Z(\delta)$ pairs $(z_{10},z_{20})$ in the case (1), or
$4$ pairs $(z_{10},z_{20})$ in the case (2). The application of the implicit function theorem
reduces to the computation of the Jacobian at $t=0$ for the system
\begin{equation}\frac{\partial\varphi_S}{\partial z_3}(\bq)=t^{-s}\frac{\partial\varphi_S}{\partial z_1}(\bq)=
t^{-s}\frac{\partial\varphi_S}{\partial z_2}(\bq)=0
\label{e-real3}\end{equation} in the case (1), or the system
$$\frac{\partial\varphi_S}{\partial z_3}(\bq)=t^{-s_1}\frac{\partial\varphi_S}{\partial z_1}(\bq)=
t^{-s_2}\frac{\partial\varphi_S}{\partial z_2}(\bq)=0$$ in the case (2). The
nondegeneracy of these Jacobians is straightforward.
\proofend

\section{Real singular surfaces in real pencils}
\label{sec-examples}
Combining the lattice path algorithm from Section \ref{sec-lp} with the
patchworking construction from Section \ref{sec-pw}, one exhibits all singular
algebraic surfaces in the given pencil, thus, making formula
(\ref{ecor03}) for the degree of the discriminant explicit. Having this in mind,
we address Problem \ref{main}(3) and give a lower bound for the maximal number of
real singular surfaces occurring in a generic real pencil in the linear system
$|{\mathcal L}_\Delta|$ (we call a pencil {\it generic} if it contains only
finitely many surfaces with singularity in $(\K^*)^3$).

\begin{theorem}\label{t-real}
For any $d\ge2$ there exists a generic real pencil of surfaces of degree $d$ in
$\PP^3$ that contains at least $\frac{3}{2}d^3+O(d^2)$ real singular surfaces.
\end{theorem}

We prove Theorem \ref{t-real} in Sections
\ref{sec-real1}--\ref{sec-rev-real}: it immediately follows from Corollaries \ref{cor-A}, \ref{cor-D}, and \ref{cor-E}. We start with defining
specific initial data for the lattice path algorithm and the patchworking construction, and then
compute the contribution of some singular tropical surfaces to the number of
real singular surfaces over the field $\K_\R$ in the corresponding pencil. By the Tarski principle, this is equivalent to the
same statement over $\R$.

\begin{remark}\label{rrev3} In principle, one can choose another directing vector for the
line through the tropical point configuration and obtain another amount of real singular surfaces
in the pencil. However, our choice has certain advantages: (1) no circuits of types {\bf B} and {\bf C}
occur, for which the patchworking construction gives a relatively small number of
real singular surfaces among all singular surfaces tropicalizing to the given tropical surface
(cf. Lemma \ref{lmt3}), (2) our choice generalizes the choice made in \cite{IKS0,IKS} for obtaining a possibly large number of real rational curves in a tropical way, (3) with our choice, one can easily enumerate almost all singular tropical surfaces
(see \cite[Appendix]{MMS2}) and compute their algebraic liftings.
\end{remark}

It follows from \cite[Corollary 6.5]{DFS}, that
$\deg\Sing(d\Delta)=4\Vol_\Z(\Delta)\cdot d^3+O(d^2)$
for any non-defective lattice polytope $\Delta$. Hence, the
lower bound of Theorem \ref{t-real} is asymptotically comparable with the
total number of (complex) singular surfaces
in the pencil.

Moreover, for an arbitrary nondegenerate convex lattice polytope $\Delta\subset\R^3$, set
\begin{equation}\alpha(\Delta)=\max\{\lambda>0\ ;\ \text{there exist}\ M\in GL(3,\Z)\ \text{and}\ v\in\R^3\ \text{such
that}\ \lambda M\Delta^3_1+v\subset\Delta\}\ .\label{e-def}\end{equation}
Here, $\Delta^3_d$ denotes the simplex in $\R^3$ with vertices $(0,0,0)$, $(d,0,0)$, $(0,d,0)$ and $(0,0,d)$.

\begin{theorem}\label{t-real1}
For an arbitrary nondegenerate convex lattice polytope $\Delta$ and any integer $d\ge1$, there exists
a generic pencil of real surfaces in
$(\C^*)^3$ with Newton polytope $d\Delta$ that
contains at least $\frac{3}{2}\alpha(\Delta)^3d^3+o(d^3)$ real singular surfaces.
\end{theorem}

This lower bound is asymptotically comparable with the degree of the discriminant too.

The proof of Theorem \ref{t-real1} is presented in Section \ref{sec-real6}. It is based on the patchworking construction in the sense
of \cite{Sh1}.

\subsection{The choice of initial data}\label{sec-real1} Denote
$$\Delta^3_d=\conv\{(0,0,0,0),(d,0,0,0),(0,d,0,0),(0,0,d,0),(0,0,0,d)\}\subset\R^3\ .$$
Fix the line $L\subset\R^3$
passing through the origin and directed by the vector $\bv=(1,\eps,\eps^2)$ with a sufficiently small
rational $\eps>0$\ \footnote{In principle, one can choose another vector
$\bv$ and find another lower bound for the number of real singular surfaces. Our choice of $\bv$ is motivated by the fact that,
in this case, possible circuits associated with lattice paths are relatively simple
(cf. \cite[Appendix A]{MMS2}).}. It then defines the following order on
$\Delta^3_d\cap\Z^3=\{\bw_k\ :\ k=0,...,N+1\}$:
\begin{equation}\bw_k=(i,j,l)\prec\bw_{k'}=(i',j',l')\quad\Longleftrightarrow\quad\begin{cases}
\text{either}\quad &i<i',\\
\text{or}\quad & i=i',\ j<j',\\
\text{or}\quad & i=i',\ j=j',\ l<l'.\end{cases}\label{eorder}\end{equation}
We shall use also the induced order
\begin{equation}(i,j)\prec(i',j')\quad\Longleftrightarrow\quad\begin{cases}\text{either}\quad & i<i',\\
\text{or}\quad & i=i',\ j<j'.\end{cases}\label{eorder1}\end{equation}
Pick $N$ points $\bx_1,...,\bx_N\in L$ satisfying (\ref{eCH}), and introduce
the configuration $\overline\bp\subset(\K^*)^3$ of $N$ points such that
\begin{equation}
\bp_i=(t^{-x_{i1}},t^{-x_{i2}},t^{-x_{i3}}),\quad \text{where}\quad \bx_i=(x_{i1},x_{i2},x_{i3}),\quad i=1,...,N
\ .\label{e-real1}\end{equation}

Notice that with these combinatorial data, no singular tropical surface with a circuit of type
{\bf B} or {\bf C} may occur (see \cite[Appendix]{MMS2}, or \cite[Section 2.D]{IKS0} for
the planar case).

\begin{lemma}\label{dop}
Let $S\in\Sing^{\tr}(\Delta^3_d,\bx)$ correspond to a lattice path $\Gamma_k$ for some $k=1,...,N+1$,
whose all segments have lattice length $1$, and let ${\mathcal S}\in\Sing(\Delta^3_d,\bp,S)$
be given by $$\varphi_S(\bz)=\sum_{i+j+l\le d}a_{ijl}z_1^iz_2^jz_3^l=0\ ,$$
where
$$a_{\bw_0}=1,\quad a_{\bw}=t^{\nu_S(\bw)}(\alpha_{\bw}+O(t^{>0})),\ \bw\in\Delta^3
\cap\Z^3\ .$$
Then
$$\alpha_{\bw_r}=\begin{cases}(-1)^r,\quad &\text{if}\ 1\le r<k,\\
(-1)^{r-1},\quad &\text{if}\ k<r\le N+1.\end{cases}$$
\end{lemma}

{\bf Proof.}
The claim immediately follows from the equations
$\varphi_S(\bp_r)=0$, $r=1,...,N$, and formulas (\ref{e-real1}).
\proofend

\subsection{Contribution of singular tropical surfaces with circuit of type {\bf A}}

\begin{lemma}\label{ld6}
Let $\bw_k=(i,d-i,0)$ with $0\le i<d$. Then, for any $5$-tuple
$$Q'_{j,l}=\{(i,d-i,0),(i,d-i-1,1),(i+1,j,0),(i+1,j-1,l),(i+1,j-1,l+1)\}\subset\Delta^3_d\ ,$$
$$j>0,\quad l\ge2,\quad j+l\le d-i-1\ ,$$
and for any $5$-tuple
$$Q''_{j,l}=\{(i,d-i,0),(i,d-i-1,1),(i+1,j,l),(i+1,j,l+1),(i+1,j-1,d-i-j)\}\subset\Delta^3_d\ ,$$
$$j>0,\quad l\ge0,\quad j+l<d-i-2\ ,$$
there exists a unique
tropical surface $S\in\Sing^\tr(\Delta^3_d,\overline\bx)$ matching the lattice path
$\Gamma_k$ (see Lemma \ref{llp1}) and having the circuit $C_S=Q'_{j,l}$, resp. $C_S=Q''_{j,l}$ of type
{\bf A}. Each of the above surfaces $S$ lifts to one real algebraic surface
${\mathcal S}\in\Sing(\Delta,\overline\bp,S)$.
\end{lemma}

{\bf Proof.}
Observe that each pentatope $\conv(Q'_{j,l})$ or $\conv(Q''_{j,l})$ is
$\Z$-affine equivalent to some $\Pi_{pq}$ defined in (\ref{epen}). Furthermore,
the last (in the sense of order (\ref{eorder})) point of $Q'_{j,l}$ is
$\bw_{n'}=(i+1,j,0)$, and the last point of $Q''_{j,l}$ is
$\bw_{n''}=(i+1,j,l+1)$, and the point $\bw_k$ is intermediate in both cases.
Furthermore, one can see that the intersection of
$\conv(Q'_{j,l})$ with $\conv\{\bw_s\ :\ 0\le s<n',\ s\ne k\}$ is a common $2$-face
$\conv\{(i,d-i-1,1),(i+1,j-1,l),(i+1,j-1,l+1)\}$, and the intersection of
$\conv(Q''_{j,l})$ with $\conv\{\bw_s\ :\ 0\le s<n'',\ s\ne k\}$ is a common
$2$-face $\conv\{(i,d-i-1,1),(i+1,j,l),(i+1,j-1,d-i-j)\}$. Furthermore, for the
points of $Q'_{j,l}$ we have the relation
$$\bw_k=(i,d-i,0)=(i+1,j,0)-l\cdot(i+1,j-1,l)
+(l-1)\cdot(i+1,j-1,l+1)+(i,d-i-1,1)\ ,$$ while for the point of $Q''_{j,l}$
$$\bw_k=(i,d-i,0)=(d-1-i-j-l)\cdot(i+1,j,l+1)-(d-2-i-j-l)\cdot(i+1,j,l)$$
$$-(i+1,j-1,d-i-j)+(i,d-i-1,1)\ ,$$ which in both cases yields, first, that $\lambda_n>0$
(in the notation of Lemma \ref{lA2}) and, second, that $|A(S,\overline\bp)|=1$
for each of the considered singular tropical surfaces $S$ (in the notation of Lemma \ref{lmt3}).
The latter relation yields that the (unique) algebraic surface
${\mathcal S}\in\Sing(\Delta^3_d,\overline\bp,S)$ is real.
\proofend

It is not difficult to show that no other surfaces $S\in\Sing^\tr(\Delta^3_d,\overline\bx)$
with a pentatopal circuit are
possible: the use of other lattice paths necessarily leads to a pair of parallel
edges in the pentatope, which is forbidden for pentatopes $\Pi_{pq}$.

\begin{corollary}\label{cor-A}
There are at least $\frac{1}{3}d^3+O(d^2)$ real algebraic surfaces ${\mathcal S}\in\Sing(\Delta^3_d,\overline\bp)$ that tropicalize
to surfaces $S\in\Sing^{\tr}(\Delta^3_d,\bx)$ with a circuit of type {\bf A}.
\end{corollary}

{\bf Proof.} It follows from Lemma \ref{ld6} that, for any triple $(i,j,l)\in\Z^3$ satisfying
$$0\le i<d,\quad 0<j,\quad 0\le l,\quad j+l\le d-2-i\ ,$$ we get a real singular algebraic surface belonging to the considered pencil. Thus, we obtain the required bound, since the number of these triples
$(i,j,l)$ is $\frac{1}{3}d^3+O(d^2)$.
\proofend

\subsection{Contribution of singular tropical surfaces with circuit of type {\bf D}}

\begin{lemma}\label{ld5}
(1) Let $\bw_k=(i,j,0)$ with $i>0$, $0<j<d-i$. Then, for any quadruple
$$Q_l=\{(i,j,0),\ (i,j,1),\ (i,j-1,l),\ (i,j-1,l+1)\}\subset\Delta^3_d,\quad
l=0,...,d-i-j,$$ there exists a unique
tropical surface $S\in\Sing^\tr(\Delta^3_d,\overline\bx)$ matching the lattice path
$\Gamma_k$ (see Lemma \ref{llp1}) and having the circuit $C_S=Q_l$ of type
{\bf D}.

(2) Let $\bw_k=(i,j,d-i-j)$ with $i>0$, $0\le j<d-i$. Then, for any quadruple
$$Q_l=\{(i,j,d-i-j),\ (i,j,d-i-j-1),\ (i,j+1,l),\ (i,j+1,l+1)\}\subset\Delta^3_d,\quad
l=0,...,d-i-j-2,$$ there exists a unique
tropical surface $S\in\Sing^\tr(\Delta^3_d,\overline\bx)$ matching the lattice path
$\Gamma_k$ and having the circuit $C_S=Q_l$ of type
{\bf D}.

(3) Let $\bw_k=(i,0,0)$ with $0<i<d$. Then, for any quadruple
$$Q_{j,l}=\{(i,0,0),\ (i,0,1),\ (i-1,j,l),\ (i-1,j,l+1)\}\subset\Delta^3_d\ ,$$
$$l=0,...,d-i-j-1, \quad j=1,\ldots,d-i\ ,$$ there exists a unique
tropical surface $S\in\Sing^\tr(\Delta^3_d,\overline\bx)$ matching the lattice path
$\Gamma_k$ and having the circuit $C_S=Q_{j,l}$ of type
{\bf D}.

(4) Each of the above surfaces $S$ satisfies $\mt(S,\overline\bx)=2$. Both singular algebraic surfaces
${\mathcal S}\in\Sing(\Delta^3_d,\bp,S)$ are real or imaginary depending on
\begin{itemize}
\item 
$\frac{1}{2}(3(d-i)+2+2j+2l)(d-i+1)\equiv1$ or $0$ modulo $2$, in case (1),
\item $\frac{1}{2}(d-i+2+2l)(d-i+1)\equiv1$ or $0$ modulo $2$, in case (2),
\item $d-i-j\equiv0$ or $1$ modulo $2$, in case (3).
\end{itemize}
\end{lemma}

{\bf Proof.}
In view of Lemmas \ref{lD1} and \ref{lD2}, to prove claims (1)-(3) one has to only show
that the quadruple $Q_l$, resp. $Q_{j,l}$ spans a parallelogram of lattice area $2$ not contained
in $\partial\Delta^3_d$,
which intersects with $\conv\{\bw_0,...,\bw_{k-1}\}$ along one of its edges, and that
the point $\bw_k$ is intermediate in $Q_l$, resp. $Q_{j,l}$ along the order
(\ref{eorder}). The relation $\mt(S,\overline\bx)=2$ follows from Lemma \ref{lmt7}. We decide on the reality of
surfaces ${\mathcal S}\in\Sing(\Delta^3_d,\bp,S)$ analysing each of the cases (1), (2), and (3) separately.

In case (1), by construction, the circuit $Q_l$ is a base of two pyramids in the dual to $S$ subdivision of $\Delta^3_d$:
one with the vertex $\bw'=(i-1,d-i+1,0)$, and the other with
the vertex $\bw''=(i+1,0,0)$. The surfaces ${\mathcal S}\in\Sing(\Delta^3_d,\bp,S)$ are real
if and only if the system (\ref{ereal1}) has two real solutions. In the considered situation, this systems
turns to be as follows.
The lattice path takes $\frac{1}{2}\cdot (j-1)\cdot (d-i+1+d-i+1-j+2)+l+1$ steps from $\bw'$ to $(i,j-1,l)$ and $\frac{1}{2}\cdot j\cdot (d-i+1+d-i+1-j+1)+2 $ steps to $(i,j,1)$.
Hence from Lemma \ref{dop}, we get (taking only the signs into account, i.e.\ reducing the number of steps mod 2 to simplify computations)
$$\alpha_{\bw'}=\sigma,\quad \alpha_{i,j-1,l}=-\alpha_{i,j-1,l+1}=
(-1)^{\lambda_1}\sigma,\quad\alpha_{ij1}=(-1)^{\lambda_2}\sigma,\quad
\alpha_{\bw''}=(-1)^{\lambda_3}\sigma\ ,$$
$$\lambda_1=l+1+(d-i)(j+1)+\frac{j^2-j}{2},\quad
\lambda_2=(d-i)j+\frac{j^2+j}{2}+1,\quad
\lambda_3=\frac{(d-i+1)(d-i+2)}{2}\ .$$
Then system (\ref{ereal1}) takes the form (without bringing the circuit to the canonical square shape)
$$z_{30}=1,\quad z_{20}=(-1)^{\lambda_2-\lambda_1},\quad
(-1)^{\lambda_3}z_{10}^2+z_{20}^{d-i+1}=0\ ,$$
and it has two solutions that are real iff $(-1)^{\lambda_3+(d-i+1)(\lambda_2-\lambda_1)}=-1$, i.e.,
$$\frac{(3(d-i)+2+2j+2l)(d-i+1)}{2}\equiv1\mod 2\ .$$

Similarly, in case (2) we again have $\bw'=(i-1,d-i+1,0)$, $\bw''=(i+1,0,0)$,
$$\alpha_{\bw'}=\sigma,\quad \alpha_{i,j,d-i-j-1}=(-1)^{\lambda_1}\sigma,\quad
\alpha_{i,j+1,l}=-\alpha_{i,j+1,l+1}=(-1)^{\lambda_2}\sigma,\quad
\alpha_{\bw''}=(-1)^{\lambda_3}\sigma\ ,$$
$$\lambda_1=(d-i)(j+1)+\frac{j^2-j}{2},\quad
\lambda_2=l+1+(d-i)(j+1)+\frac{j^2-j}{2},\quad
\lambda_3=\frac{(d-i+1)(d-i+2)}{2}\ .$$
Respectively, system (\ref{ereal1}) takes the form
$$z_{30}=1,\quad z_{20}=(-1)^{\lambda_2-\lambda_1+1},\quad
(-1)^{\lambda_3}z_{10}^2+z_{20}^{d-i+1}=0\ ,$$ and
it has two real solutions that are real iff $(-1)^{\lambda_3+(\lambda_2-\lambda_1+1)(d-i+1)}=-1$, i.e.,
$$\frac{(d-i+1)(d-i+2+2l)}{2}\equiv1\mod2\ .$$

Finally, in case (3) we get $\bw'=(i-1,j-1,d-i-j+2)$,
$\bw''=(i-1,j+1,0)$,
$$\alpha_{\bw'}=\sigma,\quad \alpha_{i-1,j,l}=-\alpha_{i-1,j,l+1}=(-1)^{\lambda_1}\sigma,\quad
\alpha_{i01}=(-1)^{\lambda_2}\sigma,\quad\alpha_{\bw''}=(-1)^{\lambda_3}\sigma\ ,$$
$$\lambda_1=l+1,\quad
\lambda_2=\frac{(d-i-j+1)(d-i-j)}{2},\quad
\lambda_3=d-i-j+2\ .$$

System (\ref{ereal1}) takes the form
$$z_{30}=1,\quad(-1)^{\lambda_2}z_{10}+(-1)^{\lambda_1+1}z_{20}^j=0,\quad
1+(-1)^{\lambda_3}z_{20}^2=0\ ,$$


and it has two solutions that are real iff $(-1)^{\lambda_3}=-1$, i.e.,
$$d-i-j+1\equiv1\mod2\ .$$
\proofend

\begin{corollary}\label{cor-D}
There are at least $\frac{1}{2}d^3+O(d^2)$ real algebraic surfaces ${\mathcal S}\in\Sing(\Delta^3_d,\overline\bp)$ that tropicalize
to surfaces $S\in\Sing^{\tr}(\Delta^3_d,\bx)$ with a circuit of type {\bf D}.
\end{corollary}

{\bf Proof.}
We claim that in each of the cases (1), (2), and (3) of
Lemma \ref{ld5}, $(1/12)d^3+O(d^2)$ singular tropical surfaces lift to pairs of real singular algebraic surfaces, and
$(1/12)d^3+O(d^2)$ singular
tropical surfaces lift to pairs of complex conjugate singular algebraic surfaces. Namely, we apply
statement (4) of Lemma \ref{ld5}. For example, in case (1), we should study the parity of the expression
$\lambda=\frac{1}{2}(3(d-i)+2+2j+2l)(d-i+1)$ in the set
$$\Lambda=\{i>0,\ 0<j<d-i,\ 0\le l\le d-i-j\}\ .$$ If $d-i$ is odd, then $\lambda\equiv\frac{1}{2}(d-i+1)$, and hence, for
$i\equiv d+1\mod 4$, we get $\lambda\equiv0\mod2$, and for $i\equiv d-1\mod4$, we get $\lambda
\equiv1\mod2$. If $d-i$ is even, then $\lambda\equiv\frac{1}{2}(3(d-i)+2+2j+2l)\mod2$, and hence, for
$j+l\equiv\frac{3}{2}(d-i)+1\mod2$, we get $\lambda\equiv0\mod2$, and for $j+l\equiv\frac{3}{2}(d-i)\mod2$, we get
$\lambda\equiv1\mod2$. Thus for (an asymptotic) half of the set $\Lambda$ we have $\lambda\equiv1\mod2$ and for the rest $\lambda\equiv0\mod2$. Thus, we obtain $\frac{1}{2}|\Lambda|+O(d^2)=(1/12)d^3+O(d^2)$ points giving singular tropical surfaces that lift to two real singular surfaces each.
The cases (2) and (3) are considered in the same manner.
\proofend

\subsection{Contribution of singular tropical surfaces with circuit of type {\bf E}}\label{sec-rev-real}

\begin{lemma}\label{ld2}
Let $\bw_k=(i,j,d-i-j)$ with $j>0$ and $i+j\le d-1$. Then there exists a unique
tropical surface $S\in\Sing^\tr(\Delta^3_d,\overline\bx)$ matching the lattice path
$\Gamma_k$ and the circuit $C_S=\{(i,j-1,d-i-j+1),(i,j,d-i-j),(i,j+1,d-i-j-1)\}$ of type
{\bf E}. Furthermore, $\mt(S,\overline\bx)=2(d-i-1)$, and all algebraic surfaces
${\mathcal S}\in\Sing(\Delta^3_d,\bp,S)$ are real.
\end{lemma}

{\bf Proof.}
The edges of the lattice path $\Gamma_k$ avoid the point $\bw_k$.
By Lemma \ref{lBCE}, the unique smooth subdivision of $\Delta^3_d$ induced by the lattice path $\Gamma_k$
defines a tropical surface $S$ with a circuit $C_S$ of type {\bf E} as indicated in the assertion.
By Lemma \ref{lmt3}(v), the set of enhancements $A(S,\bp)$ contains two elements. We claim that
they both are real. Indeed, denoting the endpoints of the circuit by
$$\bw'=(i,j-1,d-i-j+1),\quad\bw''=(i,j+1,d-i-j-1)\ ,$$ and setting $\alpha_{\bw'}=1$, we
obtain from Lemma \ref{dop}, that $\alpha_{\bw''}=(-1)^{2(d-i-j)}=1$, and hence
by formula (\ref{e-real2}) both values of $\alpha_{\bw_k}$ are real.

To allocate the singular points of the tropical surface $S$, we consider the projection
$pr^{\ba}_{x,y}:\R^3\to\R^2$ onto the $(x,y)$-plane parallel to the vector
$\ba=(0,1,-1)$. The point $\pr^{\ba}_{x,y}(C_S)=(i,d-i)$ belongs to
$\partial\Delta^2_d$, and hence the situation of Lemma \ref{lmt5}(ii) is not possible.
Set $b_0=\nu_S(i,j-1,d-i-j+1)$ and $b_1=\nu_s(i,j+1,d-i-j-1)$, and note that
\begin{equation}
0<\nu_S(i',j',l')\ll b_0\ll b_1\ll \nu_S(i'',j'',l'')\quad\text{when}\quad(i',j'+l')\prec(i,d-i)\prec(i'',j''+l'')\ .
\label{eorder4}\end{equation}
The suitably modified construction of Lemma \ref{lmt5} yields
$\Lambda'(y)=\frac{b_1-b_0}{2}(y-j)+\frac{b_0+b_1}{2}$, and that the graph of the function
$\nu_{x,y}$ is the lower convex hull of the set of points
$(\omega,-c'_\omega)$, $\omega\in\Delta^2_d\cap\Z^2\setminus\{(i,d-i)\}$, where
due to (\ref{eorder4}) we have
\begin{equation}
-c'_{i',j'}=\begin{cases}\nu_S(i',j',0)-\frac{b_1-b_0}{2}(j'-j)-\frac{b_0+b_1}{2},
\quad & \text{if}\ (i',j')\preceq(i,j)\\
\nu_S(i,j+1,j'-j-1)-b_1,\quad & \text{if}\ i'=i,\ j<j'<d-i,\\
\nu_S(i',0,j')+\frac{b_1-b_0}{2}j-\frac{b_0+b_1}{2},\quad & \text{if}\ i'>i.\end{cases}
\label{erev2}\end{equation}
One can see that all the points $(\omega,-c'_\omega)$, $\omega\in\Delta^2_d\cap\Z^2\setminus\{(i,d-i)\}$,
are vertices of the graph of $\nu_{x,y}$, and that the subdivision of
$\Delta^2_d$ induced by $\nu_{x,y}$ is a smooth triangulation built on the
lattice path, which goes through the points $\omega\in\Delta^2_d\cap\Z^2\setminus\{(i,d-i)\}$ in the order
(\ref{eorder1}). This subdivision contains the triangles $T_{j'}=\conv\{(i,d-i-1),(i+1,j'),(i+1,j'+1)\}$, $0\le j'\le d-i-2$, satisfying the conditions of
Lemma \ref{lmt5}(i). Moreover, the functions $\Lambda'':\Delta^2_d\to\R$,
linearly extending $\nu_{x,y}\big|_{T_{j'}}$, satisfy
\begin{eqnarray}\Lambda''(i,d-i)&=&-c'_{i,d-i-1}-c'_{i+1,j'+1}+c'_{i+1,j'}\nonumber\\
&=&\nu_S(i,j+1,d-i-j-2)+\nu_S(i+1,0,j'+1)-\nu_S(i+1,0,j')-b_1\nonumber\\
&=&\nu_S(i+1,0,j'+1)+o(\nu_S(i+1,0,j'+1))>0;\ \nonumber\end{eqnarray}
if $j<d-i-1 $; if $j=d-i-1$, we obtain
\begin{eqnarray}\Lambda''(i,d-i)&=&-c'_{i,d-i-1}-c'_{i+1,j'+1}+c'_{i+1,j'}\nonumber\\
&=&\nu_S(i,d-i-1,0)+\nu_S(i+1,0,j'+1)-\nu_S(i+1,0,j')-\frac{b_0+b_1}{2}\nonumber\\
&=&\nu_S(i+1,0,j'+1)+o(\nu_S(i+1,0,j'+1))>0\ .\nonumber\end{eqnarray}
Thus, by Lemma \ref{lmt5}(i), each triangle $T_{j'}$, $0\le j'\le d-i-2$, gives rise to a singular point
of the tropical surface $S$.

Hence, by Lemma \ref{lmt11}(1), we get
$\mt(S,\overline\bx)=2(d-i-1)$. Furthermore, all algebraic surfaces ${\mathcal S}\in
\Sing(\Delta^3_d,\bp,S)$ are real, since each given enhancement and a tropical singular point
as above give rise to the unique singular algebraic surface which is a real solution
of a nondegenerate real system (\ref{e-real3}).
\proofend

\begin{corollary}\label{cor-E}
There are at least $\frac{2}{3}d^3+O(d^2)$ real algebraic surfaces ${\mathcal S}\in\Sing(\Delta^3_d,\overline\bp)$ that tropicalize
to surfaces $S\in\Sing^{\tr}(\Delta^3_d,\bx)$ with a circuit of type {\bf E}.
\end{corollary}

We notice that there are $2d^3+O(d^2)$ more singular surfaces ${\mathcal S}\in\Sing(\Delta^3_d,
\bp)$ that tropicalize to singular tropical surfaces dual to subdivisions of
$\Delta^3_d$ with circuit of type {\bf E} and a lattice path containing the circuit
(see \cite[Lemmas A.1, A.2, and A.3]{MMS2}). However, the patchworking procedure
presented in Step 1 of the proof of Lemma \ref{lmt11} does not allow one to
easily decide whether the obtained singular surfaces are real.


\subsection{Proof of Theorem \ref{t-real1}}\label{sec-real6}

Let $M\in GL(3,\Z)$ and $v\in\R^3$ realize the maximal value $\lambda=\alpha(\Delta)$ in the definition
(\ref{e-def}). Since the the count of complex and real singular surfaces for a given $\Delta$ does not depend on the
$GL(3,\Z)$-action, we can suppose that $M$ is the identity.
Note that $\alpha(\Delta)\Delta^3_d+dv$ is the maximal volume simplex inscribed into $d\Delta$, and there exist
$d'=d'(d)\in \Z$ and $v'=v'(d)\in\Z^3$ such that \begin{equation}\Delta':=\Delta^3_{d'}+v'\subset\alpha(\Delta)\Delta^3_d+dv
\quad\text{and}\quad \alpha(\Delta)d-d'=O(1)\ .\label{e-real4}\end{equation}
By Theorem \ref{t-real}, there exists a configuration $\overline\bp_0\subset(\R^*)^3$ of $N_0=|\Delta^3_{d'}|-2$ points such that
the surfaces of degree $d'$ in $\PP^3$ passing through $\overline\bp_0$ form a pencil, and this pencil contains
$m=\frac{3}{2}(d')^3+O((d')^2)$ real singular surfaces. Let these surfaces be given be given by
polynomials $F^{(0)}_i\in\R[x,y,z]$ with Newton polytope $\Delta^3_{d'}+v'$, $i=1,...,m$.

Observe that by construction the above pencil intersects the discriminatal hypersurface in $\quad$
\mbox{$\deg\Sing(\Delta^3_{d'})=4(d'-1)^3$} distinct points. That is,
all the intersections are transversal, and hence by a small variation of the
configuration $\overline\bp$ we can make all truncations $(F^{(0}_i)^\delta$ of the polynomial $F^{(0)}_i$ on the faces
$\delta$ of $\Delta^3_{d'}+v'$ to be nondegenerate (i.e., defining smooth hypersurfaces in $(\C^*)^3$)
for each $i=1,...,m$).

Now, using the version of the patchworking construction from \cite[Theorem 3.1]{Sh1}, we extend the above pencil and
singular hypersurfaces $F^{(0)}_i=0$, $i=1,...,m$, to a real pencil of hypersurfaces
in the linear system $|{\mathcal L}_\Delta|$ on the
toric variety $\Tor_\K(\Delta)$ and respectively $m$ real singular hypersurfaces in it.
Since $\Vol(\alpha(\Delta)\Delta^3_d)-\Vol(\Delta^3_{d'})=O(d^2)$ (cf. (\ref{e-real4})), we then get
$m=\frac{3}{2}\alpha(\Delta)^3d^3+O(d^2)$ real singular surfaces in the pencil constructed, as required
in Theorem \ref{t-real1}.

To apply \cite[Theorem 3.1]{Sh1}, we define appropriate initial data:
\begin{enumerate}\item[(1)] {\it Combinatorial data.} Let $\delta\subset(\Delta^3_{d'}+v')$ be a two-face, $L(\delta)\subset\R^3$ the affine
plane spanned by $\delta$, $\bn_\delta\in\Z^3$ the primitive integral outer normal,
$\mu_\delta$ the value of the linear functional $x\in\R^3\mapsto\langle x,\bn_\delta\rangle\in\R$
on $\delta$. Define
$$\nu_\delta:d\Delta\to\R,\quad\nu_\delta(x)=\begin{cases}0,\quad&\text{if}\ \langle x,\bn_\delta\rangle\le\mu_\delta,\\
\langle x,\bn_\delta\rangle,\quad&\text{if}\ \langle x,\bn_\delta\rangle\ge\mu_\delta,\end{cases}$$
and set
$$\nu:d\Delta\to\R,\quad\nu=\sum_\delta\nu_\delta\ ,$$ where $\delta$ runs over all two-faces of $\Delta^3_{d'}+v'$.
Observe that $\nu$ is a convex piecewise-linear function on $d\Delta$, integral valued at $d\Delta\cap\Z^3$, and
its linearity domains divide $d\Delta$ into the union of lattice $3$-polytopes
$\Delta_0\cup...\cup\Delta_r$,
$\Delta_0=\Delta^3_{d'}+v'$. Denote by ${\mathcal G}$ the adjacency graph of the polytopes $\Delta_i$, $i=0,...,r$,
and orient ${\mathcal G}$ without oriented cycles so that $\Delta_0$ will be a pure source. This, in particular, defines
a partial order on the polytopes of the subdivision, and we will assume that the numbering
$\Delta_0,...,\Delta_r$ extends this partial order to a linear one.
\item[(2)] {\it Algebraic data.} For any $i=1,...,r$, let $N_i=|(\Delta_i\setminus\bigcup_{j<i}\Delta_j)\cap\Z^3|$ and
(if $N_i>0$) choose a generic configuration of $N_i$ points $\overline\bp_i\subset(\R^*)^3$.
Note that $$N_0+...+N_r=|d\Delta\cap\Z^3|-2=\dim\Sing(d\Delta)\ .$$
Due to the general position of each configuration $\overline\bp_i$, $i=1,...,r$, for any $j=1,...,m$, there exists
a unique sequence of polynomials $F^{(0)}_j,F^{(1)}_j,...,F^{(r)}_j\in\R[x,y,z]$. Namely, given a subsequence
$F^{(0)}_j,...,F^{(k)}_j$, $k<r$, we define $F^{(k+1)}_j$ to be the polynomial, whose coefficients of the monomials
$x^{\omega_1}y^{\omega_2}z^{\omega_3}$, $(\omega_1,\omega_2,\omega_3)\in\Delta_{k+1}\cap
\Delta_l$ coincide with the corresponding coefficients of $F^{(l)}_j$ for all $l=0,...,k$, and such that
$F^{(k+1)}_j\big|_{\overline\bp_{k+1}}\equiv0$. Due to the general position of the configurations
$\overline\bp_l$, $0\le l\le k+1$, the polynomial $F^{(k+1)}_j$ is defined uniquely, and it defines a smooth hypersurface in $(\C^*)^3$.

In addition, we define a configuration of $N=|\Delta\cap\Z^3|-2$ points in $(\K_\R^*)^3$:
$$\overline\bp=\overline\bp_0\cup\bigcup_{i=1}^r\overline\bp_i^t\ ,$$ where, for each $i=1,...,r$, the configuration
$\overline\bp_i^t$ is obtained from $\overline\bp_i\subset(\R^*)^3\subset(\K_R^*)^3$ by applying the map
$$(x,y,z)\mapsto(xt^{-\gamma_x},yt^{-\gamma_y},zt^{-\gamma_z}),\quad \nu\big|_{\Delta_i}(x,y,z)=
\gamma_xx+\gamma_yy+\gamma_zz+\gamma_0\ .$$
\item[(3)] {\it Transversality conditions.} The transversality conditions required in \cite[Theorem 3.1]{Sh1}
reduce to the following statements, which we have by construction:
\begin{itemize}
\item each polynomial $F_j^{(0)}$, $1\le j\le m$, defines a uninodal surface in $(\C^*)^3$, which
corresponds to a transverse intersection point of the
pencil defined by the configuration $\overline\bp_0$ and of the discriminant $\Sing(\Delta_0)$;
\item each polynomial $F^{(k)}_j$, $1\le j\le m$, $1\le k\le r$, is uniquely determined by the linear conditions
to have given coefficients at the points $\omega\in\Delta_k\cap\bigcup_{l<k}\Delta_l$ and to vanish at
$\overline\bp_k$.
\end{itemize}
\end{enumerate}

Thus, by \cite[Theorem 3.1]{Sh1}, each sequence $F^{(0}_j,F^{(1)}_j,...,F^{(r)}_j\in\R[x,y,z]$, $1\le j\le m$,
produces a real singular surface in the toric variety $\Tor_\K(\Delta)$ passing through the configuration
$\overline\bp\subset(\K_\R^*)^3$, which completes the proof.

\section{The asymptotically maximal number of singular points on a singular tropical surface}
\label{sec-lower}
It is well-known that the discriminant $\Sing(\Delta)$ is birationally (i.e. generically one-to-one)
covered by the incidence variety
$$\widetilde\Sing(\Delta)=\{({\mathcal S},\bq)\in|{\mathcal L}_\Delta\times\Tor_\K(\Delta)\ :\ {\mathcal S}\in\Sing(\Delta),\ \bq\in\Sing(\Sigma)\}\ .$$
It was noticed in \cite[Example 3]{MMS2} that it is not valid in the tropical setting, that is, some maximal-dimensional cones
of $\Sing^\tr(\Delta)$ are multiply covered by cones of the tropical incidence variety.
We state the problem to find a sharp upper
bound for this multiplicity.

In particular, Lemma \ref{lmt4} yields that a singular point of
a singular tropical surface with a circuit of type {\bf C} may have up to $3$ different locations, while Lemma \ref{lmt5}
does not impose any absolute upper bound to the number of such locations on singular tropical surfaces
with a circuit of type
{\bf E}. Furthermore, in the proof of Lemma \ref{ld2}, we exhibit singular tropical surfaces of
degree $d$, being of maximal-dimensional geometric type and having $d-1$ singular points.

We intend to show that one can achieve a twice better (in the asymptotic sense) bound for an
arbitrary non-degenerate Newton polytope $\Delta$. Denote by $\ld(\Delta)$ the \emph{lattice diameter} of $\Delta$, i.e. the maximal lattice length of a segment with vertices in $\Delta\cap\Z^3$. Put
$$\ld_\infty(\Delta)=\lim_{n\to\infty}\frac{\ld(n\Delta)}{n}$$
(the limit, clearly, exists and is always positive).

\begin{theorem}\label{t-lower}
For an arbitrary non-degenerate convex lattice polytope $\Delta\subset\R^3$ and any $n\ge1$ there exists
a singular tropical surface $S_n$ with Newton polytope $n\Delta$, which is of maximal-dimensional geometric type and has at least $2n\cdot\ld_\infty(\Delta)+O(1)$ tropical singular points, where $O(1)$ is a bounded function depending only on $\Delta$.
\end{theorem}

{\bf Proof.}
Our strategy is as follows. Using the construction in the proof of Lemma \ref{ld2}, we
introduce an auxiliary convex lattice polytope $\widetilde\Delta_d$ with $\ld(\widetilde\Delta_d)
=d$ and a tropical surface $\widetilde S_d\in\Sing^\tr(\widetilde\Delta_d)$ of maximal-dimensional geometric type
having $2(d-1)$ singular points. Then we inscribe $\widetilde\Delta_d$
with $d\sim n\cdot\ld_\infty(\Delta)$ into $n\Delta$ and extend the subdivision dual to $\widetilde S_d$ of $\widetilde\Delta_d$ to a subdivision of $n\Delta$ by a convex triangulation involving all lattice points as vertices, obtaining finally the required tropical surface $S_n$.

\medskip

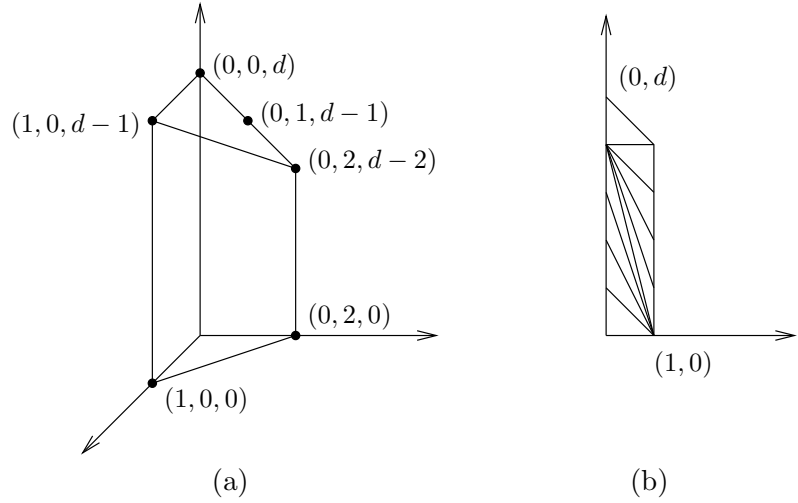
\begin{figure}
\input{figaux.pstex_t}
\caption{The polygon $\Delta'$ and the subdivision of the polygon $\delta$.}
\label{fig:aux}
\end{figure}

{\it (1)} In the hypotheses of Lemma \ref{ld2}, set $i=0$, $j=1$, $\bw_k=(0,1,d-1)$, and consider the lattice path $\Gamma_k$ that defines a subdivision of $\Delta_d^3$ with a circuit
$$C=\{(0,0,d),\ (0,1,d-1),\ (0,2,d-2)\}$$
of type {\bf E}. Introduce the subpolytope $\Delta'\subset\Delta_d^3$ given by
$$\Delta'=\conv(C\cup\{(0,0,0),(1,0,0),(0,2,0),(1,d-1,0)\}$$
(see Figure \ref{fig:aux}(a)) and define a convex piece-wise linear
function $\nu':\Delta'\to\R$ taking as its graph the lower convex hull of the points
$(\omega,\nu_S(\omega))\in\R^4$, $\omega\in\Delta'\cap\Z^3$, with $\nu_S:\Delta_d^3\to\R$ the
convex piece-wise linear function from the proof of Lemma \ref{ld2}.
It is easy to see that $\nu'$ defines a tropical surface
$S'$ with Newton polytope $\Delta'$ and a circuit $C$ (of type {\bf E}).

Observe that the surface $S'$ is of maximal-dimensional geometric type, and it has $d-1$ tropical singular points.
Indeed, the projection $\pr_{\overline a}:\R^3\to\R^2$ as in the proof of Lemma \ref{ld2} takes
$\Delta'$ onto the quadrangle $\delta=\conv\{(0,0),(1,0),(0,d),(1,d-1)\}$, on which we similarly obtain
a convex piece-wise linear function $\nu'_{x,y}$ that coincides with the restriction to $\delta$ of
the function $\nu_{x,y}:\Delta_d^2\to\R$ constructed in the proof of Lemma \ref{ld2}.
Indeed, formulas (\ref{erev2}) yield that the values of $\nu_{x,y}$ in $\delta$ are determined
merely by the values of $\nu_S$ in the integral points of $\Delta'$. Hence,
$\nu'$ defines on $\delta$ the triangulation shown in Figure \ref{fig:aux}(b), and we get the required properties of $S'$ literally in the same way as those of the surface $S$ in Lemma
\ref{ld2}.

We can correct the function $\nu':\Delta'\to\R$ by a linear function so that it will attain its
strict minimum only along the circuit $C$. We also can assume the values $\nu'(1,*,*)$ are much larger than
the values $\nu'(0,*,*)$.

Let $\rho:\R^3\to\R^3$ be the reflection on the coordinate plane $\{(0,\lambda,\mu)\ :\ \lambda,\mu\in\R\}$, and $\rho':\R^2\to\R^2$ the analogous reflection imposed on the plane by the projection $\pr_{\overline a}$ (i.e.\ the reflection on the line $\{(0,\lambda)\ :\ \lambda \in \R\}$). Define $\widetilde\Delta_d=\Delta'\cup\rho(\Delta')$ and extend the function
$\nu':\Delta'\to\R$ onto $\rho(\Delta')$ as $\rho^*\nu'$. It is clear that the resulting
function $\widetilde\nu_d:\widetilde\Delta_d\to\R$ is convex, defines on
$\rho(\Delta')$ the subdivision symmetric to that of
$\Delta'$, and furthermore, via the projection $\pr_{\overline a}$, it induces a convex piece-wise linear function
$\widetilde\nu_{x,y}:\pr_{\overline a}(\widetilde\Delta_d)=\delta\cup\rho'(\delta)\to\R$, a $\rho'$-invariant extension of $\nu'$. Thus, we obtain
a tropical surface $\widetilde S_d\in\Sing^\tr(\widetilde \Delta_d)$ of maximal-dimensional geometric type with a circuit of type
{\bf E} and $2(d-1)$ singular points.

\medskip

{\it (2)} There exists $n_0>0$, depending only on $\Delta$, such that, for any $n\ge n_0$,
the polytope $n\Delta$ suitably transformed by an automorphism of $\Z^3$ contains a subpolytope
$\widetilde\Delta_d$ with \mbox{$d=n\cdot\ld_\infty(\Delta)+O(1)$}, where $O(1)$ is a bounded function depending only on $\Delta$.
Omitting routine technicalities, we shortly comment on this claim. First, one can show that, for $n\ge n_1$ with some
fixed $n_1$, there is a vertex (that can be chosen to be the origin for all $n\Delta$) and a fixed line $L$
through the origin, on which lie the lattice segments of the maximal lattice length in $n\Delta$. Second, assuming that $L$ coincides with
the axis $\{(0,0,t),\ t\in \R\}$, one can find an integral vector $a$ such that, for any $n\ge n_2$ with some fixed $n_2$, the intersection
$(a+L)\cap n\Delta$ contains an integral segment $\sigma$ of length, which differs from the length of the maximal integral segment in $n\Delta$
by a bounded function $O(1)$ depending only on the angles between $L$ and the faces of $\Delta$, and such that the shifts of $\sigma$ by vectors $(x,y,
0)$, $|x|,|y|\le2$, lie inside $n\Delta$.
This, clearly yields the discussed claim.

\smallskip

Now, we extend the function $\widetilde\nu_d:\widetilde\Delta_d\to\R$ in a convex piece-wise linear manner to
$\nu_n:n\Delta\to\R$ as follows. We totally order the set $(n\Delta\setminus\widetilde\Delta_d)\cap\Z^3$ so that
$$\bm\not\in\conv\Big(\widetilde\Delta_d\cup\big\{\bm'\in(n\Delta\setminus\widetilde\Delta_d)\cap\Z^3\ :\ \bm'\prec\bm\big\}\Big)
\quad\text{for all}\ \bm\in(n\Delta\setminus\widetilde\Delta_d)\cap\Z^3\ ,$$ and then subsequently define the values
$\nu_n(\bm)$, $\bm\in(n\Delta\setminus\widetilde\Delta_d)\cap\Z^3$, so that
\begin{equation}
\nu_n(\bm)\gg\max\widetilde\nu_d\quad\text{and}\quad\nu_n(\bw)\gg\max\{\nu_n(\bm')\ :\ \bm'\in(n\Delta\setminus\widetilde\Delta_d)\cap\Z^3,
\ \bm'\prec\bm\}\ .\label{erev3}\end{equation}
Clearly, the complement of $\widetilde\Delta_d$ in $n\Delta$ is triangulated with vertices at all the points of $(n\Delta\setminus\widetilde\Delta_d)\cap\Z^3$,
in particular, the dual surface $S_n\in\Sing^\tr(n\Delta)$ is of maximal-dimensional geometric type with
a circuit of type {\bf E}, and it has (at least) $2(d-1)$ singular points. For the latter claim,
we observe that condition (\ref{erev3}) ensures that the convex piece-wise linear function
induced on $\pr_{\overline a}(n\Delta)$ along the rules of Lemma \ref{lmt5} is such that its restriction to
$\pr_{\overline a}(\widetilde\Delta_d)=\delta\cup\rho'(\delta)$ coincides with $\widetilde\nu_{x,y}$, and hence we obtain
$2(d-1)$ singular points as in the previous step.
\proofend

\end {document}

%% file: 3circ.pstex_t
\begin{picture}(0,0)%
\includegraphics{3circ.pstex}%
\end{picture}%
\setlength{\unitlength}{3947sp}%
\begingroup\makeatletter\ifx\SetFigFont\undefined%
\gdef\SetFigFont#1#2#3#4#5{%
  \reset@font\fontsize{#1}{#2pt}%
  \fontfamily{#3}\fontseries{#4}\fontshape{#5}%
  \selectfont}%
\fi\endgroup%
\begin{picture}(5042,992)(3064,-2966)
\put(7876,-2911){\makebox(0,0)[lb]{\smash{{\SetFigFont{10}{12.0}{\familydefault}{\mddefault}{\updefault}{\color[rgb]{0,0,0}(E)}%
}}}}
\put(3451,-2911){\makebox(0,0)[lb]{\smash{{\SetFigFont{10}{12.0}{\familydefault}{\mddefault}{\updefault}{\color[rgb]{0,0,0}(A)}%
}}}}
\put(4951,-2911){\makebox(0,0)[lb]{\smash{{\SetFigFont{10}{12.0}{\familydefault}{\mddefault}{\updefault}{\color[rgb]{0,0,0}(B)}%
}}}}
\put(6301,-2911){\makebox(0,0)[lb]{\smash{{\SetFigFont{10}{12.0}{\familydefault}{\mddefault}{\updefault}{\color[rgb]{0,0,0}(C)}%
}}}}
\put(7126,-2911){\makebox(0,0)[lb]{\smash{{\SetFigFont{10}{12.0}{\familydefault}{\mddefault}{\updefault}{\color[rgb]{0,0,0}(D)}%
}}}}
\end{picture}%

%% file: figaux.pstex_t
\begin{picture}(0,0)%
\includegraphics{figaux.pstex}%
\end{picture}%
\setlength{\unitlength}{3947sp}%
\begingroup\makeatletter\ifx\SetFigFont\undefined%
\gdef\SetFigFont#1#2#3#4#5{%
  \reset@font\fontsize{#1}{#2pt}%
  \fontfamily{#3}\fontseries{#4}\fontshape{#5}%
  \selectfont}%
\fi\endgroup%
\begin{picture}(4524,3156)(2839,-4555)
\put(6301,-4486){\makebox(0,0)[lb]{\smash{{\SetFigFont{11}{13.2}{\familydefault}{\mddefault}{\updefault}{\color[rgb]{0,0,0}(b) }%
}}}}
\put(3676,-1861){\makebox(0,0)[lb]{\smash{{\SetFigFont{10}{12.0}{\familydefault}{\mddefault}{\updefault}{\color[rgb]{0,0,0}$(0,0,d)$}%
}}}}
\put(3226,-2236){\makebox(0,0)[rb]{\smash{{\SetFigFont{10}{12.0}{\familydefault}{\mddefault}{\updefault}{\color[rgb]{0,0,0}$(1,0,d-1)$}%
}}}}
\put(3376,-3961){\makebox(0,0)[lb]{\smash{{\SetFigFont{10}{12.0}{\familydefault}{\mddefault}{\updefault}{\color[rgb]{0,0,0}$(1,0,0)$}%
}}}}
\put(4276,-3436){\makebox(0,0)[lb]{\smash{{\SetFigFont{10}{12.0}{\familydefault}{\mddefault}{\updefault}{\color[rgb]{0,0,0}$(0,2,0)$}%
}}}}
\put(3976,-2161){\makebox(0,0)[lb]{\smash{{\SetFigFont{10}{12.0}{\familydefault}{\mddefault}{\updefault}{\color[rgb]{0,0,0}$(0,1,d-1)$}%
}}}}
\put(4276,-2461){\makebox(0,0)[lb]{\smash{{\SetFigFont{10}{12.0}{\familydefault}{\mddefault}{\updefault}{\color[rgb]{0,0,0}$(0,2,d-2)$}%
}}}}
\put(6226,-1936){\makebox(0,0)[lb]{\smash{{\SetFigFont{10}{12.0}{\familydefault}{\mddefault}{\updefault}{\color[rgb]{0,0,0}$(0,d)$}%
}}}}
\put(6451,-3736){\makebox(0,0)[lb]{\smash{{\SetFigFont{10}{12.0}{\familydefault}{\mddefault}{\updefault}{\color[rgb]{0,0,0}$(1,0)$}%
}}}}
\put(3676,-4486){\makebox(0,0)[lb]{\smash{{\SetFigFont{11}{13.2}{\familydefault}{\mddefault}{\updefault}{\color[rgb]{0,0,0}(a) }%
}}}}
\end{picture}%